# SELECTING LIKELIHOOD WEIGHTS BY CROSS-VALIDATION[1]

By Xiaogang Wang and James V. Zidek

*York University and University of British Columbia*

The (relevance) weighted likelihood was introduced to formally embrace a variety of statistical procedures that trade bias for precision. Unlike its classical counterpart, the weighted likelihood combines all relevant information while inheriting many of its desirable features including good asymptotic properties. However, in order to be effective, the weights involved in its construction need to be judiciously chosen. Choosing those weights is the subject of this article in which we demonstrate the use of cross-validation. We prove the resulting weighted likelihood estimator (WLE) to be weakly consistent and asymptotically normal. An application to disease mapping data is demonstrated.

**1. Introduction.** The weighted likelihood (WL for short) has been developed for a variety of purposes. Moreover, it shares its underlying purpose with other methods such as weighted least squares and kernel smoothers which can reduce an estimator's variance while increasing its bias to reduce mean-squared error (MSE), that is, increase its precision. However, the achievement of these gains depends on choosing the weights well, which is the subject of this article. More specifically, we show that they may be data dependent (i.e., "adaptive") and chosen by *cross-validation*. The idea of data-dependent weights goes back at least to the celebrated James–Stein estimator, a WL estimator with adaptive weights that does successfully trade bias for variance [Hu and Zidek (2002)].

To describe the WL, we assume independent random response vectors $\mathbf{X}_1, \ldots, \mathbf{X}_m$ with probability density functions $f_1(\cdot;\theta_1), \ldots, f_m(\cdot;\theta_m)$, where $\mathbf{X}_i = (X_{i1}, \ldots, X_{in_i})^t$. Further suppose that only population 1, in particular

Received June 2003; revised May 2004.
[1]Supported in part by the Natural Sciences and Engineering Research Council of Canada.

$\theta_1$, an unknown vector of parameters, is of inferential interest. Given data $\mathbf{X}=\mathbf{x}$, the classical likelihood would be

$$L_1(\mathbf{x_1},\theta_1)=\prod_{j=1}^{n_1} f(x_{1j};\theta_1).$$

When the remaining parameters $\theta_2,\dots,\theta_m$ are thought to resemble $\theta_1$, the WL is defined as

$$\mathrm{WL}(\mathbf{x};\theta_1)=\prod_{i=1}^{m}\prod_{j=1}^{n_i} f_1(x_{ij};\theta_1)^{\lambda_i},$$

where $\boldsymbol{\lambda}=(\lambda_1,\dots,\lambda_m)$, the "weights vector," must be specified. Notice that the parameters from the remaining populations, $\theta_2,\dots,\theta_m$, unlike the data they generate, do not appear in the WL, since inferential interest focuses on $\theta_1$. It follows that

$$\log \mathrm{WL}(\mathbf{x};\theta_1)=\sum_{i=1}^{m}\sum_{j=1}^{n_i}\lambda_i \log f_1(x_{ij};\theta_1).$$

The WL extends the local likelihood method of Tibshirani and Hastie (1987) for nonparametric regression, although the idea predates them [see Hu and Zidek (2002) for a review]. Following Hu (1997), Hu and Zidek (1995, 2001, 2002) extend the local likelihood to a more general setting. However, the aim is the same. Their method also combines all relevant information in samples from populations thought to resemble the one of inferential interest.

The maximum WL estimator (WLE) for $\theta_1$, say $\tilde{\theta}_1$, is defined by

$$\tilde{\theta}_1=\arg\sup_{\theta_1\in\Theta}\mathrm{WL}(\mathbf{x};\theta_1).$$

In many cases the WLE may be obtained by solving the *estimating equation*:

$$(\partial/\partial\theta_1)\log\mathrm{WL}(\mathbf{x};\theta_1)=0.$$

Note that uniqueness of the WLE is not assumed.

Like the MLE, the WLE has a number of desirable properties [Hu and Zidek (2002)], in particular consistency and asymptotic normality under reasonable general conditions [Hu (1997) and Wang, van Eeden and Zidek (2004)]. However, these asymptotic properties have only been shown with fixed weights and hence need to be extended in this article to cover the estimators we obtain using cross-validation.

In its most primitive but nevertheless useful form, the cross-validation procedure consists of controlled and uncontrolled division of the data sample into two subsamples. For example, a subsample can be generated by deleting one or more observations or it can be a random sample from the data set. Stone (1974) began the systematic study of cross-validatory choice and assessment in statistical prediction. Both Stone (1974) and Geisser (1975) discuss its application to the *K-group* problem and use a linear combination of the sample means from different groups to estimate a common mean. Breiman and Friedman (1997) also demonstrate the benefit of using cross-validation to obtain linear combinations of predictors that perform well in multivariate regression.

The article is organized as follows. The adaptive weights are derived in Section 2. The asymptotic properties of the resulting WLE are presented in Section 3. Results of simulation studies are discussed in Section 4. In Section 5 an application to disease mapping data demonstrates the benefits of using the proposed method in conjunction with the WLE when compared with traditional estimators.

**2. Choosing adaptive weights.** For cross-validation there are many ways of dividing the entire sample into subsets, such as a random selection technique. However, we use the simplest *leave-one-out* approach in this article since the analytic forms of the optimum weights are then completely tractable for the linear WLE. Denote the vector of parameters and the weight vector by $\boldsymbol{\theta} = (\theta_1, \theta_2, \dots, \theta_m, \rho)$ and $\boldsymbol{\lambda} = (\lambda_1, \lambda_2, \dots, \lambda_m)$, respectively. Let $\boldsymbol{\lambda}_e^{\text{opt}}$ and $\boldsymbol{\lambda}_u^{\text{opt}}$ denote the optimum weight vectors to be defined in the sequel for samples with equal and unequal sizes, respectively. We require that $\sum_{i=1}^m \lambda_i = 1$ in this section and throughout this article.

Suppose that we have $m$ populations which might be related to each other. The probability density functions or probability mass functions are of the form $f_i(x; \theta_i)$ with $\theta_i$ as the parameter vector for population $i$. Assume that

$$\begin{array}{cccccc} X_{11}, & X_{12}, & X_{13}, & \dots, & X_{1n_1} & \overset{\text{i.i.d.}}{\sim} f_1(x;\theta_1), \\ X_{21}, & X_{22}, & X_{23}, & \dots, & X_{2n_2} & \overset{\text{i.i.d.}}{\sim} f_2(x;\theta_2), \\ \vdots & & & & & \vdots \\ X_{m1}, & X_{m2}, & X_{m3}, & \dots, & X_{mn_m} & \overset{\text{i.i.d.}}{\sim} f_m(x;\theta_m), \end{array}$$

where, for fixed $i$, the $\{X_{ij}\}$ are observations obtained from population $i$ and so on. Assume that observations obtained within each population are independent and identically distributed. Also observations from one population are independent of those from other populations except that $\text{Corr}(X_{ij}, X_{kj}) = \rho$, for any fixed $j$ and $i \neq k$. That is, observations having the same second subscripts are not necessarily independent even though they are from different populations. This would allow a spatial correlation structure but not a temporal one. We also assume that $E(X_{ij}) = \phi(\theta_i) = \phi_i$, say, for $j = 1, 2, \dots, n_i, i = 1, 2, \dots, m$. The population parameter of the first population, $\theta_1$, is of inferential interest.

Our cross-validatory approach of estimating the weights for the WLE flows from taking prediction as our inferential objective. In other words, we seek an estimator $\hat{\theta}_1$ of $\theta_1$ that enables us to predict accurately, in some sense, a randomly drawn element $X_1^*$ from the first population. But how should the precision of $\hat{\theta}_1$ be assessed?

One answer is the expected log score. Denoting by "$E$" the expectation with respect to the conditional distribution of $X_1$ given $\theta_1$, that score is $E[\log f_1(X_1|\hat{\theta}_1)]$, an index of $\hat{\theta}_1$'s performance.

However, the complexity of that index makes its use impractical in applications such as that in Section 5. We therefore adopt an approximation as a compromise. To obtain the approximation, we assume a one-to-one mapping of $\theta_1$ into $(\phi_1, \tau_1)$ where the range of $\phi_1$ covers that of $X_1$. In fact, with an abuse of notion we represent $\theta_1$ by $\theta_1 = (\phi_1, \tau_1)$ and $\hat{\theta}_1$ in a similar way. We further assume that

$$\left.\frac{\partial \log E[f_1(X_1|\hat{\theta}_1)]}{\partial \hat{\phi}_1}\right|_{\hat{\theta}_1=\theta_1} = 0$$

and

$$\left.\frac{\partial^2 \log E[f_1(X_1|\hat{\theta}_1)]}{\partial^2 \hat{\phi}_1}\right|_{\hat{\theta}_1=\theta_1} < 0,$$

for all $\theta$ with all higher-order derivatives being assumed to exist. These assumptions are satisfied for the normal distribution, for example, and more importantly for our application in Section 5, the Poisson distribution.

Under these assumptions, the first-order term in a three-term Taylor expansion of the expected log score vanishes. Therefore, ignoring irrelevant terms and factors, we obtain $(\hat{\phi}_1 - \phi_1)^2$ as an approximation to the negative expected log score as a measure of $\hat{\phi}_1$'s precision. Finally, for its empirical assessment, we estimate the unknown $\phi_1$ in this measure by $X_1$. Moreover, we adopt that empirical measure to obtain adaptive weights by cross-validation. To that end, we use $(-j)$ to indicate that the $j$th item has been dropped from the sample.

Taking the usual path, we predict $X_{1j}$ by $\phi(\tilde{\theta}_1^{(-j)})$, the WLE of its mean without using the $X_{1j}$. Note that $\phi(\tilde{\theta}_1^{(-j)})$ is a function of the weight vector $\boldsymbol{\lambda}$ by the construction of the WLE. Based on the log score approximation above, a natural measure for the discrepancy of the WLE becomes

$$D(\boldsymbol{\lambda}) = \sum_{j=1}^{n_1} (X_{1j} - \phi(\tilde{\theta}_1^{(-j)}))^2. \tag{1}$$

The optimum weights are derived such that the minimum of $D(\boldsymbol{\lambda})$ is achieved for fixed sample sizes $n_1, n_2, \ldots, n_m$ and $\sum_{i=1}^m \lambda_i = 1$.

If the inferential interest is on the means of some commonly used distributions from the exponential family, such as normal, binomial, exponential and Poisson distributions, it then follows that $\phi(\tilde{\theta})$ is simply a linear combination of those MLEs for each population. In this section we will investigate the behavior of the optimum weights by cross-validation for the linear case since we can derive the analytical forms of the optimum weights from (1).

2.1. *Linear WLEs for equal sample sizes.* Stone (1974) and Geisser (1975) discuss the application of the cross-validation approach to the so-called *K-group* problem. Suppose that the data set $S$ consists of $n$ observations in each of $K$ groups. The mean predictor for the $i$th group is

$$\hat{\mu}_i = (1-\alpha)\overline{X}_{i\cdot} + \alpha\overline{X}_{\cdot\cdot},$$

where $\overline{X}_{i\cdot} = \frac{1}{n}\sum_{j=1}^{n} X_{ij}$ and $\overline{X}_{\cdot\cdot} = \frac{1}{K}\sum_{i=1}^{m}\overline{X}_{i\cdot}$. If our interest focuses on group 1, the relevant predictor is

$$\hat{\mu}_1 = \left(1 - \frac{K-1}{K}\alpha\right)\overline{X}_{1\cdot} + \sum_{i=2}^{m}\frac{\alpha}{K}\overline{X}_{i\cdot},$$

where $\alpha$ is a parameter. Stone (1974) uses cross-validation to derive an optimal value for $\alpha$. We remark that the above formula is just a particular linear combination of the sample means.

We consider more general linear combinations and throughout this section assume $n_1 = n_2 = \cdots = n_m = n$. Let $\tilde{\theta}_1^{(e)}$ denote the WLE obtained through cross-validation. If $\phi(\theta) = \theta$, the linear WLE for $\theta_1$ is then defined as

$$\tilde{\theta}_1^{(e)} = \sum_{i=1}^{m}\lambda_i\overline{X}_{i\cdot},$$

where $\sum_{i=1}^{m}\lambda_i = 1$.

In this section we will use cross-validation by simultaneously deleting $X_{1j}, X_{2j}, \ldots, X_{mj}$ for each fixed $j$. That is, we delete one data point from each sample at each step. This might be appropriate if these data points are obtained at the same time point and strong associations exist among them. By simultaneously deleting $X_{1j}, X_{2j}, \ldots, X_{mj}$ for each fixed $j$, we might achieve numerical stability of the cross-validation procedure. An alternative approach is to delete a data point from only the first sample at each step. That approach will be studied in this section as well.

Let $\overline{X}_{i\cdot}^{(-j)}$ be the sample mean of the $i$th sample with $j$th element in that sample excluded. A natural measure for the discrepancy of $\tilde{\theta}_1$ might be

$$\begin{aligned} D_e^{(m)} &= \sum_{j=1}^{n}\left(X_{1j} - \sum_{i=1}^{m}\lambda_i\overline{X}_{i\cdot}^{(-j)}\right)^2 \\ &= c(\underline{\mathbf{X}}) - 2\boldsymbol{\lambda}^t b_e(\underline{\mathbf{X}}) + \boldsymbol{\lambda}^t A_e(\underline{\mathbf{X}})\boldsymbol{\lambda}, \end{aligned}$$

where $c(\underline{\mathbf{X}}) = \sum_{j=1}^n X_{1j}^2$, $(b_e(\underline{\mathbf{X}}))_i = \sum_{j=1}^n X_{1j}\overline{X}_{i\cdot}^{(-j)}$ and $(A_e(\underline{X}))_{ik} = \sum_{j=1}^n \overline{X}_{i\cdot}^{(-j)}\overline{X}_{k\cdot}^{(-j)}, i = 1,2,\dots,n, k = 1,2,\dots,m$. For expository simplicity, let $b_e = b_e(\underline{X})$ and $A_e = A_e(\underline{X})$ in this article.

An optimum weight vector by the cross-validation procedure is defined to be a vector that minimizes the objective function $D_e^{(m)}$ and satisfies $\sum_{i=1}^m \lambda_i = 1$.

2.1.1. *Two-population case.* For simplicity, first consider the simple case of just two populations,

$$
\begin{aligned}
&X_{11}, \quad X_{12}, \quad X_{13}, \quad \dots, \quad X_{1n} \overset{\text{i.i.d.}}{\sim} f_1(x;\theta_1,\sigma_1^2),\\
&X_{21}, \quad X_{22}, \quad X_{23}, \quad \dots, \quad X_{2n} \overset{\text{i.i.d.}}{\sim} f_2(x;\theta_2,\sigma_2^2),
\end{aligned}
$$

with $E(X_{1j}) = \theta_1$, $E(X_{2j}) = \theta_2$, $\mathrm{Var}(X_{1j}) = \sigma_1^2$ and $\mathrm{Var}(X_{2j}) = \sigma_2^2$. Furthermore, assume that $\rho = \mathrm{cor}(X_{1j}, X_{2j})$, $j = 1,2,\dots,n$. Denote $\boldsymbol{\theta}^0 = (\theta_1^0, \theta_2^0)$ where $\theta_1^0$ and $\theta_2^0$ are the true values for $\theta_1$ and $\theta_2$, respectively.

We seek the optimum weights $\lambda_1$ and $\lambda_2$ with $\lambda_1 + \lambda_2 = 1$ such that they minimize the following objective function:

$$
D_e^{(2)} = \sum_{j=1}^n (X_{1j} - \lambda_1\overline{X}_{1\cdot}^{(-j)} - \lambda_2\overline{X}_{2\cdot}^{(-j)})^2 - \gamma(\lambda_1 + \lambda_2 - 1).
$$

Differentiating $D_e^{(2)}$ with respect to $\lambda_1$ and $\lambda_2$, we have

$$
\begin{aligned}
\frac{\partial D_e^{(2)}}{\partial \lambda_1} &= -\sum_{j=1}^n \overline{X}_{1\cdot}^{(-j)}(X_{1j} - \lambda_1\overline{X}_{1\cdot}^{(-j)} - \lambda_2\overline{X}_{2\cdot}^{(-j)}) - \gamma = 0,\\
\frac{\partial D_e^{(2)}}{\partial \lambda_2} &= -\sum_{j=1}^n \overline{X}_{2\cdot}^{(-j)}(X_{1j} - \lambda_1\overline{X}_{1\cdot}^{(-j)} - \lambda_2\overline{X}_{2\cdot}^{(-j)}) - \gamma = 0.
\end{aligned} \tag{2}
$$

It follows that

$$
\begin{aligned}
\lambda_1^{\mathrm{opt}}(\mathbf{X}) &= 1 - \frac{\sum_{j=1}^n(\overline{X}_{1\cdot}^{(-j)} - \overline{X}_{2\cdot}^{(-j)})(\overline{X}_{1\cdot}^{(-j)} - X_{1j})}{\sum_{j=1}^n(\overline{X}_{1\cdot}^{(-j)} - \overline{X}_{2\cdot}^{(-j)})^2},\\
\lambda_2^{\mathrm{opt}}(\mathbf{X}) &= 1 - \lambda_1^{\mathrm{opt}}(\mathbf{X}).
\end{aligned} \tag{3}
$$

LEMMA 2.1. *The following identity holds*:

$$
\lambda_1^{\mathrm{opt}} = 1 - \lambda_2^{\mathrm{opt}} \quad \textit{and} \quad \lambda_2^{\mathrm{opt}} = S_2^e / S_1^e,
$$

*where*

$$
\begin{aligned}
S_1^e &= \frac{n(n-2)}{(n-1)^2}(\overline{X}_{1\cdot} - \overline{X}_{2\cdot})^2 + \frac{1}{n(n-1)^2}\sum_{j=1}^n (X_{1j} - X_{2j})^2,\\
S_2^e &= \frac{n}{(n-1)^2}(\hat{\sigma}_1^2 - \widehat{\mathrm{cov}}),
\end{aligned}
$$

*where* $\hat{\sigma_1}^2 = \frac{1}{n}\sum_{j=1}^n (X_{1j} - \overline{X}_{1\cdot})^2$ *and* $\widehat{\text{cov}} = \frac{1}{n}\sum_{j=1}^n (X_{1j} - \overline{X}_{1\cdot})(X_{2j} - \overline{X}_{2\cdot})$.

The value of $\lambda_2^{\text{opt}}$ can be seen as some sort of measure of relevance between the two samples. If this "measure" is almost zero, the formula for $\lambda_2^{\text{opt}}$ will assume a very small value. This implies that there is no need to combine the two populations if the difference between the two sample means is relatively large or the second sample has little relevance to the first one. The weights chosen by the cross-validation procedure will then guard against the undesirable scenario in which too much bias might be introduced into the estimation procedure. On the other hand, if the second sample does contain valuable information about the parameter of interest, then the cross-validation procedure will recognize that by assigning a nonzero value to $\lambda_2^{\text{opt}}$. Note that knowledge of the variances and correlation is not assumed.

PROPOSITION 2.1. *If* $\rho < \frac{\sigma_1}{\sigma_2}$, *then*

$$P_{\boldsymbol{\theta}^0}(\lambda_2^{\text{opt}} > 0) \xrightarrow{P_{\boldsymbol{\theta}^0}} 1.$$

We remark that the condition $\rho < \sigma_1/\sigma_2$ is satisfied if $\sigma_2 < \sigma_1$ or $\rho < 0$. If the condition $\rho < \sigma_1/\sigma_2$ is not satisfied, then $\lambda_2^{\text{opt}}$ will have a negative sign for sufficiently large $n$. However, the value of $\lambda_2^{\text{opt}}$ will converge to zero as shown in the next proposition.

PROPOSITION 2.2. *If* $\theta_1^0 \neq \theta_2^0$, *then, for any given* $\varepsilon > 0$,

$$P_{\boldsymbol{\theta}^0}(|\lambda_1^{\text{opt}} - 1| \le \varepsilon) \longrightarrow 1 \quad \textit{and} \quad P_{\boldsymbol{\theta}^0}(|\lambda_2^{\text{opt}}| < \varepsilon) \longrightarrow 1.$$

The asymptotic limit of the weights will not exist if $\theta_1^0$ equals $\theta_2^0$. The cross-validation procedure will not be able to detect the difference of the two populations if there is none. This problem can be solved by defining $\lambda_2^{\text{opt}} = \frac{S_2^e}{S_1^e + \delta_e}$ where $\delta_e$ is a small positive constant.

2.1.2. *Alternative matrix representation of the optimum weights.* In order to handle more than two populations, it is necessary to derive an alternative matrix representation of $\lambda^{\text{opt}}$. Define $e_n = \frac{n}{n-1}$. It can be verified that

$$\begin{aligned}\overline{x}_{i\cdot}^{(-j)}\overline{x}_{k\cdot}^{(-j)} &= \left(e_n\overline{x}_{i\cdot} - \frac{1}{n-1}x_{ij}\right)\left(e_n\overline{x}_{k\cdot} - \frac{1}{n-1}x_{kj}\right) \\ &= e_n^2\overline{x}_{i\cdot}\overline{x}_{k\cdot} - \frac{e_n}{n-1}x_{ij}\overline{x}_{k\cdot} - \frac{e_n}{n-1}x_{kj}\overline{x}_{i\cdot} + \left(\frac{1}{n-1}\right)^2 x_{ij}x_{kj}.\end{aligned}$$

Thus, we have

$$\sum_{j=1}^{n} \overline{x}_{i\cdot}^{(-j)} \overline{x}_{k\cdot}^{(-j)} = \left(e_n^2(n-2) + \frac{e_n}{n-1}\right)\hat{\theta}_i\hat{\theta}_k + \frac{e_n}{n-1}\widehat{\mathrm{cov}}_{ik}^2, \tag{4}$$

where

$$\hat{\theta}_i = \overline{x}_{i\cdot}, \qquad i = 1, 2, \ldots, m,$$
$$\widehat{\mathrm{cov}}_{ik} = \frac{1}{n}\sum_{j=1}^{n}(x_{ij} - \overline{x}_{i\cdot})(x_{kj} - \overline{x}_{k\cdot}).$$

Recall that, for $1 \le i \le m$ and $1 \le k \le m$,

$$A_{e(ik)} = \sum_{j=1}^{n} \overline{x}_{i\cdot}^{(-j)} \overline{x}_{k\cdot}^{(-j)}.$$

It follows that

$$A_e = \frac{e_n}{n-1}\hat{\Sigma} + \left(e_n^2(n-2) + \frac{e_n}{n-1}\right)\hat{\boldsymbol{\theta}}\hat{\boldsymbol{\theta}}^t, \tag{5}$$

where $\Sigma_{ik} = \widehat{\mathrm{cov}}_{ik}$ and $\hat{\boldsymbol{\theta}} = (\overline{x}_{1\cdot}, \ldots, \overline{x}_{m\cdot})$.

We also have

$$b_{e(i)}(\mathbf{x}) = A_{1i} - \frac{e_n}{n-1}\sum_{j=1}^{n}(x_{1j} - \overline{x}_{1\cdot})x_{ij}. \tag{6}$$

It then follows that

$$b_e(\mathbf{x}) = A_1 - e_n^2\widehat{\Sigma}_1, \tag{7}$$

where $A_1$ is the first column of $A_e$ and $\widehat{\Sigma}_1$ is the first column of the sample covariance matrix $\widehat{\Sigma}$. We are now in a position to derive the optimum weights in matrix form when the sample sizes are equal.

PROPOSITION 2.3. *The optimum weight vector which minimizes $D_e^{(m)}$ takes the form*

$$\boldsymbol{\lambda}_e^{\mathrm{opt}} = (1, 0, 0, \ldots, 0)^t - e_n^2\left(A_e^{-1}\widehat{\Sigma}_1 - \frac{\mathbf{1}^t A_e^{-1}\widehat{\Sigma}_1}{\mathbf{1}^t A_e^{-1}\mathbf{1}}A_e^{-1}\mathbf{1}\right).$$

We remark that $A_e$ is invertible since $\hat{\Sigma}$ is invertible. Note that the expression of the weight vector in the two-population case can also be derived by using the matrix representation given as above.

2.2. *Linear WLE for unequal sample sizes.* In the previous section we discussed choosing the optimum weights when the samples sizes are equal. In this section we propose to use the cross-validation procedure for choosing adaptive weights for unequal sample sizes. If the sample sizes are not equal, it is not clear whether the *delete-one-column* approach is reasonable. For example, suppose that there are 10 observations in the first sample and there are 5 observations in the second. Then there is no observation to delete for the second sample for the second half of the cross-validation steps. Furthermore, we might lose accuracy in prediction by deleting one entire column when sample sizes are small. Thus, we propose an alternative method that deletes only one data point from the first sample and keeps all the data points from the rest of the samples when the sample sizes are not equal.

2.2.1. *Two-population case.* Let us again consider the two populations. The optimum weights $\boldsymbol{\lambda}_u^{\text{opt}}$ are obtained by minimizing the objective function

$$D_u^{(2)}(\boldsymbol{\lambda}) = \sum_{j=1}^{n_1}(X_{1j} - \lambda_1 \overline{X}_{1\cdot}^{(-j)} - \lambda_2 \overline{X}_{2\cdot})^2,$$

where $\sum_{i=1}^m \lambda_i = 1$ and $\overline{X}_{1\cdot}^{(-j)} = \frac{1}{n_1-1}\sum_{k\neq j}^{n_1} X_{1k}$. We remark that the major difference between $D_e^{(2)}$ and $D_u^{(2)}$ is that only the $j$th data point of the first sample is left out for the $j$th term in $D_u^{(2)}$.

Under the condition that $\lambda_1 + \lambda_2 = 1$, we can rewrite $D_u^{(2)}$ as a function of $\lambda_1$:

$$\begin{aligned} D_u^{(2)} &= \sum_{j=1}^{n_1}(X_{1j} - \lambda_1 \overline{X}_{1\cdot}^{(-j)} - (1-\lambda_1)\overline{X}_{2\cdot})^2 \\ &= \sum_{j=1}^{n_1}((X_{1j} - \overline{X}_{2\cdot}) + \lambda_1(\overline{X}_{2\cdot} - \overline{X}_{1\cdot}^{(-j)}))^2. \end{aligned}$$

By differentiating $D_u^{(2)}$ with respect to $\lambda_1$, we then have

$$\lambda_1^{\text{opt}} = \frac{n_1(\overline{X}_{1\cdot} - \overline{X}_{2\cdot})^2 - (n_1/(n_1-1))\hat{\sigma}_1^2}{n_1(\overline{X}_{1\cdot} - \overline{X}_{2\cdot})^2 + (n_1/(n_1-1)^2)\hat{\sigma}_1^2}, \qquad \lambda_2^{\text{opt}} = 1 - \lambda_1^{\text{opt}}. \tag{8}$$

The adaptive optimum weights still converge to $(1, 0)$ when the sample sizes are not equal.

PROPOSITION 2.4. *If* $\theta_1^0 \neq \theta_2^0$, *then* $\lambda_1^{\text{opt}} \xrightarrow{P_{\boldsymbol{\theta}^0}} 1$ *and* $\lambda_2^{\text{opt}} \xrightarrow{P_{\boldsymbol{\theta}^0}} 0$.

2.2.2. *Optimum weights by cross-validation.* We now derive the matrix representation for the optimum weights by cross-validation when the sample sizes are not all equal. The objective function is defined as follows:

$$\begin{aligned} D_u^{(m)} &= \sum_{j=1}^{n_1}\left(X_{1j} - \lambda_1 \overline{X}_{1\cdot}^{(-j)} - \sum_{i=2}^{m} \lambda_i \overline{X}_{i\cdot}\right)^2 \\ &= c(\underline{X}) - 2\mathbf{b}(\underline{X})\boldsymbol{\lambda}_u + \boldsymbol{\lambda}_u^t A(\underline{X})\boldsymbol{\lambda}_u, \end{aligned}$$

where

$$\begin{aligned} b_1 &= \sum_{j=1}^{n_1} X_{1j}\left(\overline{X}_{1\cdot} + \frac{1}{n_1-1}(\overline{X}_{1\cdot} - X_{1j})\right) = n_1 \overline{X}_{1\cdot}^2 - \frac{n_1}{n_1-1}\hat{\sigma}_1^2, \\ b_i &= n_1 \overline{X}_{1\cdot}\overline{X}_{i\cdot}, \qquad i = 2, \dots, m, \end{aligned}$$

and

$$\begin{aligned} a_{11} &= \sum_{j=1}^{n_1}\left(\overline{X}_{1\cdot} + \frac{1}{n_1-1}(\overline{X}_{1\cdot} - X_{1j})\right)^2 = n_1 \overline{X}_{1\cdot}^2 + \frac{n_1}{(n_1-1)^2}\hat{\sigma}_1^2, \\ a_{ij} &= n_1 \overline{X}_{i\cdot}\overline{X}_{j\cdot}, \qquad i \neq 1 \text{ or } j \neq 1. \end{aligned}$$

It then follows that

$$A = n_1(\hat{\theta}_1, \hat{\theta}_2, \dots, \hat{\theta}_m)^t(\hat{\theta}_1, \hat{\theta}_2, \dots, \hat{\theta}_m) + D,$$

where

$$\begin{aligned} d_{11} &= \frac{n_1}{(n_1-1)^2}\hat{\sigma}_1^2, \\ d_{ij} &= 0, \qquad i \neq 1 \text{ or } j \neq 1. \end{aligned}$$

By the elementary rank inequality, it follows that

$$\operatorname{rank}(A) \le \operatorname{rank}(\hat{\theta}^t\hat{\theta}) + \operatorname{rank}(D) = 2.$$

It implies that

$$\operatorname{rank}(A) < m \qquad \text{if } m > 2.$$

Since $A$ is not invertible for $m > 2$, the Lagrange method will not work in this case. The $g$-inverse of the matrix $A$ could be used instead.

## 3. Asymptotic properties of the weights.

In this section we present the asymptotic properties of the cross-validated weights for the general case. Let $\hat{\theta}_1^{(n_1)}$ be the MLE based on the first sample of size $n_1$. Let $\hat{\theta}_1^{(-j)}$ and $\tilde{\theta}_1^{(-j)}$ be the MLE and WLE, respectively, based on $m$ samples without the $j$th data point from the first sample. This generalizes the two cases where either only the $j$th data point is deleted from the first sample or the $j$th data point from each sample is deleted. Note that $\tilde{\theta}_1^{(-j)}$ is a function of the weights function $\boldsymbol{\lambda}$. Let $\frac{1}{n_1}D_{n_1}$ be the average discrepancy in the cross-validation given by

$$\frac{1}{n_1}D_{n_1}(\boldsymbol{\lambda}) = \frac{1}{n_1}\sum_{j=1}^{n_1}(X_{1j} - \phi(\tilde{\theta}_1^{(-j)}))^2.$$

Let $\boldsymbol{\lambda}^{(cv)}$ be the optimum weights chosen by cross-validation. We require that $\sum_{i=1}^{m}\lambda_i = 1$. Let $\boldsymbol{\theta}^0 = (\theta_1^0, \theta_2, \dots, \theta_m)$, where $\theta_1^0$ is the true value of $\theta_1$. We then have the following theorem.

THEOREM 3.1. *Assume that*:

(i) $\frac{1}{n_1}D_{n_1}$ *has a unique minimum for any fixed* $n_1$;

(ii) $\frac{1}{n_1}\sum_{j=1}^{n_1}(\phi(\hat{\theta}_1^{(-j)}) - \phi(\theta_1^0)) \xrightarrow{P_{\boldsymbol{\theta}^0}} 0$ *as* $n_1 \to \infty$;

(iii) $P_{\boldsymbol{\theta}^0}(\frac{1}{n_1}\sum_{j=1}^{n_1}(X_{1j} - \phi(\hat{\theta}_1^{(-j)}))^2 < K) \xrightarrow{P_{\boldsymbol{\theta}^0}} 1$ *for some constant* $0 < K < \infty$;

(iv) $P_{\boldsymbol{\theta}^0}(|\phi(\hat{\theta}_1^{n_1}) - \phi(\tilde{\theta}_1^{n_1})| > M) = o(\frac{1}{n_1})$ *for some constant* $0 < M < \infty$. *Then*

$$\boldsymbol{\lambda}^{(cv)} \xrightarrow{P_{\boldsymbol{\theta}^0}} \mathbf{w}_0 = (1, 0, 0, \dots, 0)^t. \tag{9}$$

The assumptions of the above theorem are satisfied by the linear-WLE case presented in Section 2. We state that fact in the following corollary whose proof is straightforward and omitted for brevity.

COROLLARY 3.1. *Suppose* $X_{i1}, X_{i2}, \dots, X_{in}$ *are independent with density function* $f(x, \theta_i)$, $i = 1, 2$. *If the WLE has linear form and* $\mu_1 \neq \mu_2$, *then*

$$\boldsymbol{\lambda}^{(cv)} \xrightarrow{P_{\boldsymbol{\theta}^0}} \mathbf{w}_0 = (1, 0)^t. \tag{10}$$

Furthermore, Theorem 3.1 also applies to cases in which the WLE does not have the linear form. One such important case involves the log-normal distribution, which is widely used in practice. Suppose $X_{ij} \overset{\text{ind.}}{\sim} \mathrm{LN}(\mu_i, 1), j = 1, \dots, n, i = 1, 2$, where $\mu_i$ and 1 denote, respectively, the mean and standard deviation of the $\log X_{ij}$ for all $i$ and $j$. It can be verified that, for $i = 1, 2$,

$$E_{\mu_1^0}(X_{ij}) = \phi(\mu_i^0) = e^{\mu_i^0 + 1/2}, \qquad j = 1, 2, \dots, n.$$

It also follows that the MLE and the WLE are given by

$$\mathrm{MLE}(\mu_1) = \hat{\mu}_1 = \frac{1}{n}\sum_{j=1}^{n}\log(x_{1j}), \tag{11}$$

$$\mathrm{WLE}(\mu_1) = \tilde{\mu}_1 = \frac{\lambda_1}{n}\sum_{j=1}^{n}\log(x_{1j}) + \frac{\lambda_2}{n}\sum_{j=1}^{n}\log(x_{2j}), \tag{12}$$

where $\lambda_1 + \lambda_2 = 1$.

Therefore,

$$\phi(\hat{\mu}_1^{(-j)}) = \exp\left\{\frac{1}{n-1}\sum_{k\neq j}\log(X_{1k}) + 1/2\right\}, \tag{13}$$

$$\phi(\tilde{\mu}_1^{(-j)}) = \exp\left\{\frac{\lambda_1}{n-1}\sum_{k\neq j}\log(X_{1k}) + \frac{\lambda_2}{n-1}\sum_{k\neq j}\log(X_{2k}) + 1/2\right\}, \tag{14}$$

for $j = 1, 2, \ldots, n$.

Therefore the average discrepancy of cross-validation for the log-normal case is given by

$$\begin{aligned}\frac{1}{n}D_n(\lambda_1,\lambda_2) = \frac{1}{n}\sum_{j=1}^{n}\Bigg(X_{1j} - \exp\Bigg\{&\frac{\lambda_1}{n-1}\sum_{k\neq j}\log(x_{1k}) \\ &+ \frac{\lambda_2}{n-1}\sum_{k\neq j}\log(x_{2k}) + 1/2\Bigg\}\Bigg)^2.\end{aligned} \tag{15}$$

Since we require that $\lambda_1 + \lambda_2 = 1$, we can rewrite the average discrepancy as

$$\frac{1}{n}D_n(1-\lambda_2,\lambda_2) = \frac{1}{n}\sum_{j=1}^{n}(X_{1j} - e^{\overline{Y}_{1\cdot}^{(-j)} + \lambda_2(\overline{Y}_{2\cdot}^{(-j)} - \overline{Y}_{1\cdot}^{(-j)}) + 1/2})^2, \tag{16}$$

where

$$\overline{Y}_{i\cdot}^{(-j)} = \frac{1}{n-1}\sum_{k\neq j}Y_{ik} \quad \text{and} \quad Y_{ij} = \log(X_{ij}), \qquad i = 1, 2, j = 1, 2, \ldots, n.$$

We then have the following lemma and corollary.

LEMMA 3.1. *Assume that $X_{i1}, X_{i2}, \ldots, X_{in}$ are independent random variables and follow the log-normal distribution with parameters $(\mu_i, 1)$, $i = 1, 2$. Let $\lambda_2^*(n)$ be the optimum weight that minimizes $\frac{1}{n}D_n(1-\lambda_2,\lambda_2)$ for any fixed $n$. If $\mu_1 \neq \mu_2$, then* (i) *$\frac{1}{n}D_n(1-\lambda_2,\lambda_2)$ is strictly convex;* (ii) *$\lim_{n\to\infty}\lambda_2^*(n)$ exists and $|\lim_{n\to\infty}\lambda_2^*(n)| < 1$ with probability* 1.

COROLLARY 3.2. *Under the assumptions of Lemma 3.1, if $\mu_1 \neq \mu_2$, then*

$$\boldsymbol{\lambda}^{(cv)} \xrightarrow{P_{\mu^0}} \mathbf{w}_0 = (1, 0)^t. \tag{17}$$

Wang, van Eeden and Zidek (2004) establish the asymptotic normality of the WLE for fixed weights. Under certain regularity conditions and by Theorem 3.1, we then have the following asymptotic results for using adaptive weights.

THEOREM 3.2. *For each* $\theta_1^0$, *the true value of* $\theta_1$, *and each* $\theta_1 \neq \theta_1^0$,

$$\lim_{n_1\to\infty} P_{\boldsymbol{\theta}^0}\left(\prod_{i=1}^{m}\prod_{j=1}^{n_i} f(X_{ij};\theta_1^0)^{\lambda_i^{(\mathbf{n})}(\mathbf{X})} > \prod_{i=1}^{m}\prod_{j=1}^{n_i} f(X_{ij};\theta_1)^{\lambda_i^{(\mathbf{n})}(\mathbf{X})}\right) = 1,$$

*for any* $\theta_2, \theta_3, \ldots, \theta_m, \theta_i \in \Theta, i = 2, 3, \ldots, m$.

THEOREM 3.3. *For any sequence of maximum weighted likelihood estimates* $\tilde{\theta}_1^{(n_1)}$ *of* $\theta_1$ *constructed with adaptive weights* $\lambda_i^{(\mathbf{n})}(\mathbf{X})$, *and for all* $\varepsilon > 0$,

$$\lim_{n_1\to\infty} P_{\boldsymbol{\theta}^0}(\|\tilde{\theta}_1^{(n_1)} - \theta_1^0\| > \varepsilon) = 0,$$

*for any* $\theta_2, \theta_3, \ldots, \theta_m, \theta_i \in \Theta, i = 2, 3, \ldots, m$.

We assume that the parameter space is an open subset of $R^p$. The asymptotic normality of the WLE constructed by cross-validated weights follows.

THEOREM 3.4 (Multidimensional). *Suppose*:

(i) *for almost all* $x$ *the first and second partial derivatives of* $f_1(x;\theta)$ *with respect to* $\theta$ *exist, are continuous in* $\theta \in \Theta$, *and may be passed through the integral sign in* $\int f_1(x;\theta)\,d\nu(x) = 1$;

(ii) *there exist three functions* $G_1(x)$, $G_2(x)$ *and* $G_3(x)$ *such that for all* $\theta_2, \ldots, \theta_m$, $E_{\theta^0}|G_l(X_{ij})|^2 \leq K_l < \infty, l = 1, 2, 3,\ i = 1, \ldots, m$, *and in some neighborhood of* $\theta_1^0$ *each component of* $\psi(x) = \frac{\partial}{\partial\theta} f_1(x;\theta)$ [*resp.* $\dot{\psi}(x)$] *are bounded in absolute value by* $G_1(x)$ [*resp.* $G_2(x)$] *uniformly in* $\theta_1 \in \Theta$. *Further,*

$$\frac{\partial^3 \log f_1(x;\theta_1)}{\partial\theta_{1k_1}\,\partial\theta_{1k_2}\,\partial\theta_{1k_3}},$$

$k_1, k_2, k_3 = 1, \ldots, p$, *are bounded by* $G_3(x)$ *uniformly in* $\theta_1 \in \Theta$;

(iii) $I(\theta_1^0)$ *is positive definite.*

*Then there exists a sequence of roots of the weighted likelihood function based on adaptive weights* $\tilde{\theta}_1^{(n_1)}$ *that is weakly consistent and*

$$\sqrt{n_1}(\tilde{\theta}_1^{(n_1)} - \theta_1^0) \xrightarrow{D} N(0, I(\theta_1^0)) \qquad \text{as } n_1 \to \infty.$$

**4. Simulation studies.** To demonstrate and verify the benefits of using cross-validation procedures described in previous sections, we perform simulations according to the following algorithm that deletes the $j$th point from each sample, that is, a *delete-one-column* approach. Let $\mu_1^0$ and $\mu_2^0$ denote the true values of the parameters. Let $C = \mu_1^0 - \mu_2^0$, which is the difference between the two population means.

TABLE 1
*MSE $*$ 100 of the MLE and the WLE and their standard deviations $*$ 100 for samples of equal sizes generated from $N(0,1)$ and $N(0.3,1)$. A correction term is employed in the calculations of the optimum weights to handle numerical instability*

| $n$ | MSE(MLE) | SD of $(\mathrm{MLE}-\theta_1^0)^2$ | MSE(WLE) | SD of $(\mathrm{WLE}-\theta_1^0)^2$ | $\frac{\mathrm{MSE(WLE)}}{\mathrm{MSE(MLE)}}$ |
|---|---|---|---|---|---|
| 10 | 10 | 15 | 8 | 12 | 80 |
| 20 | 4 | 6 | 4 | 5 | 85 |
| 30 | 3 | 4 | 3 | 4 | 87 |
| 40 | 3 | 4 | 2 | 3 | 91 |
| 50 | 2 | 3 | 2 | 2 | 92 |
| 60 | 2 | 2 | 2 | 2 | 94 |

STEP 1. Draw random samples of size $n$ from $f_1(x;\mu_1^0)$ and $f_2(x;\mu_2^0)$.

STEP 2. Calculate the cross-validated optimum weights by using (3).

STEP 3. Calculate $(\mathrm{MLE}-\mu_1^0)^2$ and $(\mathrm{WLE}-\mu_1^0)^2$.

Repeat Steps 1–3, 1000 times. Calculate the averages and standard deviations of the squared estimation error differences for both the MLE and WLE. Calculate the averages and standard deviations of the optimum weights.

We generate random samples from $N(\mu_1^0,\sigma_1^2)$ and $N(\mu_2^0,\sigma_2^2)$ where we set $\sigma_1=\sigma_2=1$ for simplicity. For the purpose of the demonstration, we set $\mu_1^0=0$ and $\mu_2^0=0.3$, which is 30% of the variance. Table 1 shows some results for the case $\mu_1^0=0$ and $\mu_2^0=0.3$. Setting $\mu_1^0=0$, we tried other values for $C$. In general, the larger the value of $C$, the less improvement in the MSE. For example, if we set $\sigma_1^0=\sigma_2^0=1$ and $C=\mu_2^0-\mu_1^0=1$, the ratio of the MSE for MLE and WLE will be almost 1. This implies that the cross-validation procedure will not make much use of the second sample in this situation.

It is obvious from Table 1 that the MSE of the WLE is much smaller than that of the MLE for small and moderate sample sizes. The standard deviations of the squared differences for the WLE are less than or equal to those of the MLE. This suggests that not only the WLE achieves smaller MSE but also its MSE has less variation than that of the MLE. Intuitively, as the sample size increases, the importance of the second sample diminishes. As indicated by Table 2, the cross-validation procedure realizes this and then assigns a larger value to $\lambda_1$ as the first sample size increases. The optimum weights do increase towards the asymptotic weights $(1,0)$ for the normal case, albeit quite slowly.

We repeat the procedure for Poisson distributions with $\mathcal{P}(3)$ and $\mathcal{P}(3.6)$. Some of the results are shown in Tables 3 and 4. The results for the Poisson distributions differ from the normal case. The most striking difference is in the ratio of the WLE's average MSE versus that of the MLE. The WLE achieves a smaller average MSE when the sample sizes are less than 30. These results contrast with the normal case, where the critical sample size is 45.

We remark that the reduction in MSE will disappear if we set $C = \mu_2^0 - \mu_1^0 = 1.5$ in the above case. Thus, the cross-validation procedure will not combine the two samples if the second sample does not help to predict the behavior of the first. We should emphasize that the value $C$ in both cases is not used in the cross-validation procedure itself.

We remark that simulations using the *delete-one-point* approach have also been done. They give quite similar results.

**5. Application to disease mapping.** In this section we address the problem of analyzing disease mapping data. In particular, we demonstrate a weighted likelihood alternative to the hierarchical Bayes approach that has been used in references cited in the discussion section. Our approach allows the data themselves to select the weights through cross-validation. We thereby avoid the (need of a prior for modeling) in order to guess the latent patterns of environmental hazards that may lead to the adverse health effects being mapped. Such hazards include air pollution that has been associated with respiratory morbidity [see, e.g., Burnett and Krewski (1994) and Zidek, White and Le (1998)].

TABLE 2
*Average optimum weights* $*\,100$ *and their standard deviations* $*\,100$ *for samples of equal sizes generated from* $N(0,1)$ *and* $N(0.3,1)$. *A correction term is employed in the calculations of the optimum weights to handle numerical instability*

| $n$ | **AVE. of $\lambda_1$** | **AVE. of $\lambda_2$** | **SD of $\lambda_1$ and $\lambda_2$** |
|---|---|---|---|
| 10 | 79 | 21 | 6 |
| 20 | 85 | 15 | 4 |
| 30 | 88 | 11 | 3 |
| 40 | 90 | 10 | 3 |
| 50 | 91 | 9 | 2 |
| 60 | 92 | 8 | 2 |

TABLE 3
*MSE* $*\,100$ *of the MLE and the WLE and their standard deviations* $*\,100$ *for samples of equal sizes generated from* $\mathcal{P}(3)$ *and* $\mathcal{P}(3.6)$. *A correction term is employed in the calculations of the optimum weights to handle numerical instability*

| $n$ | **MSE(MLE)** | **SD of $(\mathrm{MLE} - \theta_1^0)^2$** | **MSE(WLE)** | **SD of $(\mathrm{WLE} - \theta_1^0)^2$** | $\frac{\mathrm{MSE(WLE)}}{\mathrm{MSE(MLE)}}$ |
|---|---|---|---|---|---|
| 10 | 31 | 45 | 27 | 40 | 86 |
| 20 | 15 | 22 | 14 | 19 | 90 |
| 30 | 10 | 14 | 9 | 13 | 94 |
| 40 | 8 | 11 | 8 | 10 | 96 |
| 50 | 6 | 8 | 5 | 8 | 97 |
| 60 | 5 | 8 | 5 | 7 | 97 |

Our demonstration involves parallel time series of weekly hospital admissions for respiratory disease in residents of 733 census subdivisions (CSD) in southern Ontario. The data are collected from the May-to-August periods from 1983 to 1988. In this demonstration we confine attention to certain densely populated areas.

Let us consider the problem of estimating the rate of weekly hospital admissions of CSD 380, the one with the largest total annual hospital admissions among all CSDs from 1983 to 1988. This proves to be a challenging task due to the sparseness of the data set. The original data set contains many 0's, representing no hospital admissions. For example, although CSD 380 has the largest total number of hospital admissions among all the CSDs, no patient was admitted during 112 out of the 123 days in the summer of 1983. On some days, however, quite a number of people sought treatment for acute respiratory disease possibly due to high levels of air pollution in their regions. Again referring to CSD 380, 17 patients were admitted on day 51 alone in 1983.

A more graphical description of these irregularities in admission counts for this CSD is seen in Figure 1. There daily counts are shown and the problems of data sparseness and high level of variations are extreme. In fact, in this demonstration we have chosen to avoid the complexities of modeling these daily data series and we turn instead to weekly counts. While those problems

Table 4

*Average optimum weights* $*\,100$ *and their standard deviations* $*\,100$ *for samples of equal sizes generated from* $\mathcal{P}(3)$ *and* $\mathcal{P}(3.6)$*. A correction term is employed in the calculations of the optimum weights to handle numerical instability* remain, they are not nearly so acute. In total, each of the summers in the years covered by our study has 17 weeks. For simplicity, the data obtained in the last few days of each summer are dropped from the analysis since they do not constitute a whole week.

| $n$ | AVE. of $\lambda_1$ | AVE. of $\lambda_2$ | SD of $\lambda_1$ and $\lambda_2$ |
|---|---|---|---|
| 10 | 80 | 20 | 7 |
| 20 | 86 | 14 | 5 |
| 30 | 88 | 12 | 4 |
| 40 | 90 | 10 | 3 |
| 50 | 92 | 8 | 3 |
| 60 | 92 | 8 | 2 |

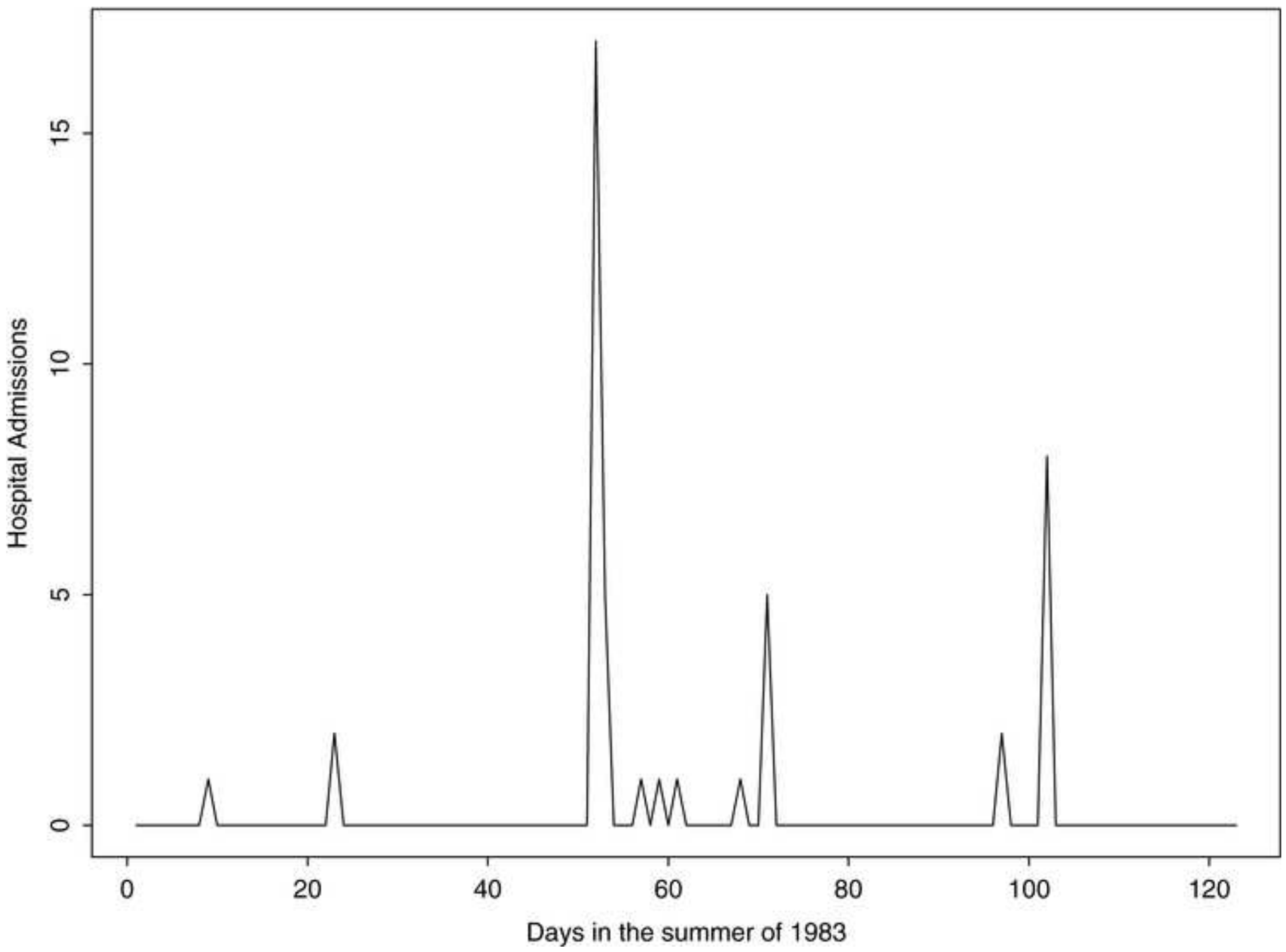

Fig. 1. *Daily hospital admissions for CSD* 380 *in the summer of* 1983.

5.1. *Weighted likelihood estimation.* We assume the weekly hospital admissions for any given CSD follow Poisson distributions, that is, for year $q$, CSD $i$ and week $j$,

$$Y_{ij}^q \overset{\text{ind.}}{\sim} \mathcal{P}(\theta_{ij}^q), \qquad j = 1, 2, \ldots, 17; i = 1, 2, \ldots, 733; q = 1, 2, \ldots, 6.$$

The raw estimates of $\theta_{ij}^q$, namely $Y_{ij}^q$, are highly unreliable since the effective sample size in this case is 1. Moreover, each CSD may contain only a small group of people who suffer respiratory diseases. These considerations point to the need to "borrow strength," a standard tool of disease mapping techniques. That is, the information in neighboring CSDs can be combined to produce more reliable estimates while introducing only a small amount of bias.

For any given CSD, the "neighboring" CSDs are defined to be CSDs in close proximity to CSD 380. To estimate the rate of weekly hospital admissions in a particular CSD, we would expect that neighboring subdivisions contain relevant information which might help us to derive a better estimate than the traditional sample average. Thus, the Euclidean distances between the target CSD and other CSDs in the data set are calculated by using the longitudes and latitudes. We apply a somewhat arbitrary threshold, 0.2, to the Euclidean distances in order to define neighbors. For CSD 380, neighboring CSDs turn out to be CSDs 362, 366 and 367.

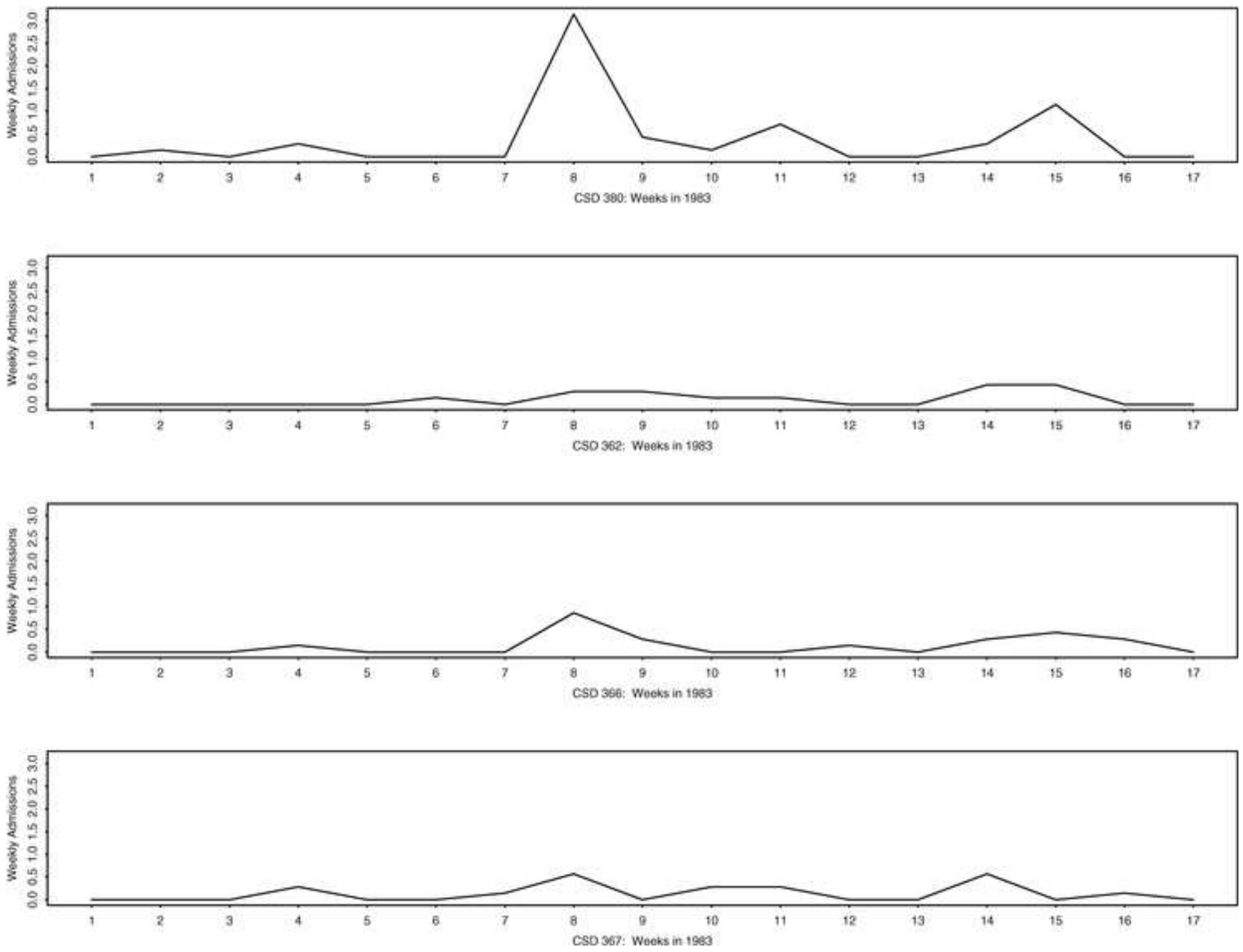

FIG. 2. *Hospital admissions for CSDs* 380*,* 362*,* 366 *and* 367 *in* 1983*.*

The time series plots of weekly hospital admissions for those selected CSDs in 1983 are shown in Figure 2. Hospital admissions of these CSDs indeed seem to be related since the major peaks in the time series plot occurred at roughly the same time points. However, as noted earlier, the data from other CSDs may introduce bias. Thus the WLE's weights are needed to control the degree of bias.

To find cross-validatory choices for these weights, we consider purely as a working assumption that $\theta^q_{ij} = \theta^q_i$ for $j = 1, 2, \ldots, 17$. In fact, that assumption does not seem tenable since every year week 8 has markedly larger numbers of hospital admissions for CSD 380 than the remaining weeks. For example, in 1983, there are 21 admissions in week 8 while the second largest weekly count is only 7 in week 15. Thus, we are forced to drop week 8 from our working assumption and instead assume $\theta^q_{ij} = \theta^q_i$ for $j = 1, 2, \ldots, 7, 9, \ldots, 17$.

In fact, the sample means and variances of the weekly hospital admissions for those 16 weeks of CSD 380 are quite close to each other, in support of our assumption. One alternative to assuming the constancy of weights over the whole summer would be the use of a moving window just a few weeks in width. We leave that option for future work.

For Poisson distributions the MLE of $\theta_1^q$ is the sample average of the weekly admissions of CSD 380, while the WLE is a linear combination of the sample averages for each CSD. Thus, the *weighted likelihood estimate* of the population mean of weekly hospital admissions for a CSD is

$$\mathrm{WLE}^q = \sum_{i=1}^{4} \lambda_i^q \overline{Y}_{i\cdot}^q, \qquad q = 1, 2, \ldots, 6,$$

where $\overline{Y}_{i\cdot}^q$ is the overall sample average of CSD $i$ for year $q$.

In our analysis the weights are selected by the cross-validation procedure proposed in Section 2. Recall that the cross-validated weights for equal sample sizes are

$$\boldsymbol{\lambda}^q = A_q^{-1}\left(b_q + \frac{1 - \mathbf{1}^t A_q^{-1} b}{\mathbf{1}^t A_q^{-1}\mathbf{1}} A_q^{-1}\mathbf{1}\right),$$

where $b_q(\underline{y}) = \sum_{j=1}^{17} Y_{1j}^q \overline{Y}_{i\cdot}^{q(-j)}$ and $A_q(\underline{y})_{ik} = \sum_{j=1}^{17} \overline{Y}_{i\cdot}^{q(-j)} \overline{Y}_{k\cdot}^{q(-j)}$, $i = 1, 2, 3, 4$; $k = 1, 2, 3, 4$.

5.2. *Results of the analysis.* We assess the performance of the MLE and the WLE by comparing their MSEs. The MSEs of the MLE and the WLE are defined by, for $q = 1, 2, \ldots, 6$,

$$\mathrm{MSE}_M^q(\theta_1^q) = E_{\theta_1^q}(\overline{Y}_{1\cdot}^q - \theta_1^q)^2,$$

$$\mathrm{MSE}_W^q(\theta_1^q) = E_{\theta_1^q}\left(\sum_{i=1}^{4} \lambda_i^q \overline{Y}_{i\cdot}^q - \theta_1^q\right)^2.$$

In fact, the $\theta_1^q$'s are unknown. We then estimate the $\mathrm{MSE}_M$ and $\mathrm{MSE}_W$ by replacing $\theta_1^q$ by the MLE. Under the assumption of Poisson distributions, the estimated MSE for the MLE is given by

$$\mathrm{MSE}_M^q = \widehat{\mathrm{var}(Y_{11})}/16, \qquad q = 1, 2, \ldots, 6.$$

The estimated MSE for the WLE is given as follows:

$$\begin{aligned}\mathrm{MSE}_W^q &= E\left(\sum_{i=1}^{m} \lambda_i^q \overline{Y}_{i\cdot}^q - \theta_1^q\right)^2 \\ &= \mathrm{Var}\left(\sum_{i=1}^{m} \lambda_i^q \overline{Y}_{i\cdot}^q\right) + \left(E\sum_{i=1}^{m} \lambda_i^q \overline{Y}_{i\cdot}^q - \theta_1^q\right)^2\end{aligned}$$

Table 5

*Estimated MSEs for the MLE and the WLE. All entries have been multiplied by 100. The MSEs have also been multiplied by 16 since there are 16 weeks*

| **Year** | **MLE** | **WLE** | $\mathbf{16 * \widehat{MSE}_M^q}$ | $\mathbf{16 * \widehat{MSE}_W^q}$ | $\mathbf{\widehat{MSE}_W^q / \widehat{MSE}_M^q}$ |
|---|---|---|---|---|---|
| 1 | 19 | 17 | 10 | 8.4 | 80 |
| 2 | 33 | 28 | 24 | 13 | 87 |
| 3 | 23 | 26 | 29 | 14 | 54 |
| 4 | 15 | 22 | 16 | 8.4 | 96 |
| 5 | 30 | 32 | 30 | 13 | 80 |
| 6 | 38 | 41 | 41 | 24 | 54 |

$$\approx \sum_{i=1}^{4}\sum_{k=1}^{4} \lambda_i^q \lambda_j^q \widehat{\text{cov}}(\overline{Y}_{i\cdot}^q, \overline{Y}_{j\cdot}^q) + \left(\sum_{i=1}^{m} \lambda_i^q \overline{Y}_{i\cdot}^q - \overline{Y}_{1\cdot}^q\right)^2.$$

The estimated MSEs for the MLE and the WLE are given in Table 5. It can be seen that the MSE for the WLE is much smaller than that of the MLE. In fact, the average reduction of the MSE by using WLE is about 25%.

Combining information across these CSDs might also help us in predictions since the patterns exhibited in one neighboring location in a particular year might manifest themselves at the location of interest the next year. To assess the performance of the WLE, we also use the WLE derived from one particular year to predict the overall weekly average of the next year. The overall prediction error is defined as the average of those prediction errors. To be more specific, the overall prediction errors for the WLE and the MLE are defined as follows:

$$\text{PRED}_M = \sqrt{\frac{1}{5}\sum_{q=1}^{5}(\overline{Y}_{1\cdot}^q - \overline{Y}_{1\cdot}^{q+1})^2},$$

$$\text{PRED}_W = \sqrt{\frac{1}{5}\sum_{q=1}^{5}(\text{WLE}^q - \overline{Y}_{1\cdot}^{q+1})^2}.$$

The average prediction error for the MLE, $\text{Pred}_M$, is 0.065, while $\text{Pred}_W$, the average prediction error for the WLE, is 0.047, which is about 72% of that of the MLE.

From Table 6, we see that there is strong linear association between CSD 380 and CSD 366. However, the weight assigned to CSD 366 is the smallest one. It shows that CSDs with higher correlations contain less information for the prediction since they might have patterns too similar to the target CSD for a given year to be helpful in the prediction for the next year. Thus CSD 366, which has the smallest weight, should not be included in the analysis. Therefore, the "neighborhood" of CSD 380 in the analysis should only include CSD 362 and CSD 367.

In general, we might examine those CSDs that are in close proximity to the target CSD. We can calculate the weight for each selected CSD by using the cross-validation procedure. CSDs with small weights should be dropped from the analysis since they are not deemed to be helpful.

The predictive distributions for the weekly totals will be Poisson as well. We can then derive the 95% predictive intervals for the weekly average hospital admissions. This might be criticized as failing to take into account the uncertainty of the unknown parameter. Smith (1999) argues that the traditional plug-in method has a small MSE compared to the posterior mean under certain circumstances. In particular, it has a smaller MSE when the true value of the parameter is not large. Let $\mathrm{CI}_W$ and $\mathrm{CI}_M$ be the 95% predictive intervals of the weekly averages calculated from the WLE and the MLE, respectively. The results are shown in Table 7.

The weighted likelihood framework discussed in this article requires the observations obtained from each population to follow the same distribution. However, including the week 8 data would violate that assumption. Including them in the analysis would have negative impact on the analysis by invalidating the homogeneity assumption of our model. Nevertheless, we re-did the analysis to see that impact. The adaptive weights and the correlation matrix for 1986 are shown in Table 8. We observe that the weight for the population of interest is almost 0. This is not acceptable since the inference will ignore the data from the first population. In this case, week 8 for CSD 380 has an observation that is almost 20 times larger than the rest of them. Since the cross-validation procedure is based on the predictive mechanism, thus it is difficult for the procedure to rely on the data points from the first population for accurate predictions. As a result, it will assign large weights to the other CSDs, especially those less correlated with the target one or having a smaller variance. Consequently, the weights will not be able to control the bias as they are designed to. Instead, they will introduce large bias into the inference.

TABLE 6
*Correlation matrix* $*\,100$ *and the weights* $*\,100$ *for* 1984

| | **CSD 380** | **CSD 362** | **CSD 366** | **CSD 367** | **Weights** |
|---|---|---|---|---|---|
| CSD 380 | 100 | 42 | 91 | 57 | 46 |
| CSD 362 | 42 | 100 | 40 | 63 | 20 |
| CSD 366 | 91 | 40 | 100 | 55 | 12 |
| CSD 367 | 57 | 63 | 55 | 100 | 22 |

TABLE 7
*Predictive confidence intervals of the MLE and the WLE for CSD* 380

| **Year** | $\mathbf{CI}_M$ | $\mathbf{CI}_W$ |
|---|---|---|
| 1983 | [0, 3] | [0, 3] |
| 1984 | [0, 5] | [0, 4] |
| 1985 | [0, 4] | [0, 4] |
| 1986 | [0, 3] | [0, 4] |
| 1987 | [0, 4] | [0, 5] |
| 1988 | [0, 5] | [0, 6] |

Table 9 presents the results obtained when the data from week 8 are dropped for 1986. As in Table 6, a large weight, about 50%, is put back onto CSD 380, the population of interest. Therefore, data from week 8 must be dropped from the analysis in order to control the bias. We discuss some alternative methods for detecting unusual weeks in the discussion section. In principle, we could fit a separate model for that week. But here it would not be feasible because of the rather small sample size. We note that the MLE and WLE are both unstable for small sample sizes although the WLE will have better performance as shown in the simulation study.

**6. Discussion and future work.** The asymptotic results established in this article are based on the assumption that the sample size of the popula-tion of interest goes to infinity. They do not apply to the situation when the sample size for the population of interest remains small or moderate while the sample sizes of other populations go to infinity. If the sample size of the population of interest is very small, say either 1 or 0, and the number of populations goes to infinity, then the asymptotic paradigm proposed by Hu (1997) would be appropriate.

TABLE 8
*Correlation matrix* $*\,100$ *and the weights* $*\,100$ *for* 1986 *when week* 8 *is included in the analysis*

| | **CSD 380** | **CSD 362** | **CSD 366** | **CSD 367** | **Weights** |
|---|---|---|---|---|---|
| CSD 380 | 100 | 88 | 74 | 22 | 0.1 |
| CSD 362 | 88 | 100 | 76 | 32 | 28 |
| CSD 366 | 74 | 76 | 100 | 44 | 30 |
| CSD 367 | 22 | 32 | 44 | 100 | 42 |

TABLE 9
*Correlation matrix* $*\,100$ *and the weights* $*\,100$ *for* 1986 *when week* 8 *is excluded in the analysis*

| | **CSD 380** | **CSD 362** | **CSD 366** | **CSD 367** | **Weights** |
|---|---|---|---|---|---|
| CSD 380 | 100 | 23 | 19 | 7.6 | 48 |
| CSD 362 | 23 | 100 | 38 | 29 | 18 |
| CSD 366 | 19 | 38 | 100 | 44 | 31 |
| CSD 367 | 7.6 | 29 | 44 | 100 | 2.6 |

There are other choices of weights function proposed in the literature. In the context of local likelihood discussed by Copas (1995), Tibshirani and Hastie (1987) and Eguchi and Copas (1998), the weight function there is essentially a kernel function with center $t$ and bandwidth $h$. Hunsberger (1994) proposes a weight function that assigns zero weight to an observation if it is outside a certain neighborhood. Since a kernel-type weight function uses Euclidean distance, it might not reflect the underlying spatial structure well as we have seen in the disease mapping example. Hu and Rosenberger (2000) propose weight functions in analyzing adaptive designs when time trends are present. They investigate two classes of weight functions, namely the exponential and polynomial types. But the weight function proposed in this article does not assume any specific functional form or rely on the choice of distance function. The adaptive weights chosen by cross-validation are data dependent and determined solely by minimizing the proposed predictive discrepancy measure.

The analysis presented in Section 5 is merely a demonstration of the weighted likelihood method. Through exploratory analysis, we find that data from week 8 are quite different from the rest of the weeks. Therefore they were dropped from the analysis. Given the high dimensionality and actual sizes of current data sets in disease mapping, it is not always practical to detect those unusual weeks by manual exploratory analysis. One automatic approach to detect patterns for the weekly data is to partition those weeks into homogeneous subgroups by using some clustering algorithms. Unlike the standard clustering in disease mapping that is normally done on the spatial grid, the grouping in our case should be done on the temporal scale. We applied a standard *K-means* algorithm with two clusters to the data set. The *K-means* clustering algorithm successfully identified week 8 as the only member of one cluster and the rest of the weeks were assigned to another cluster. When the number of clusters is unknown, it then must be estimated. The estimation of number of clusters is a very difficult problem in cluster analysis. It is beyond the scope of this article. Fraley and Raftery (1998) discuss the problem of determining the structure of clustered data without prior knowledge of the number of clusters. Cheeseman and Stutz (1996) propose an algorithm, the so-called *AutoClass*, that can estimate the number of clusters and then perform the partition. Once the partition is achieved, the weighted likelihood method can then be applied to those clusters separately. One of our future works is how to combine the results from different clusters in a sensible way. Furthermore, the spatial structure is incorporated into the weighted likelihood through the adaptive weights. However, the current model cannot handle temporal structures. One natural extension of the proposed weighted likelihood framework is to extend it to handle both spatial and temporal structures.

Bayes methods including empirical and hierarchical Bayes methods are widely used in analyzing disease mapping data. Manton et al. (1989) discuss the empirical Bayes procedures for stabilizing maps of cancer mortality rates. Ghosh, Natarajan, Waller and Kim (1999) propose a very general hierarchical Bayes spatial generalized model that is considered broad enough to cover a large number of situations where spatial structures need to be incorporated. In particular, they propose the following:

$$\theta_i = q_i = x_i^t \mathbf{b} + u_i + v_i, \qquad i = 1, 2, \ldots, m,$$

where the $q_i$ are known constants, $x_i$ are covariates, $u_i$ and $v_i$ are mutually independent with $v_i \overset{\text{i.i.d.}}{\sim} N(0, \sigma_v^2)$ and the $u_i$ have joint probability density function

$$f(u) \propto ((\sigma_u)^2)^{-1/2m} \exp\left( -\sum_{i=1}^{m} \sum_{j \neq i} (u_i - u_j)^2 w_{ij} / (2\sigma_u^2) \right),$$

where $w_{ij} \geq 0$ for all $1 \leq i \neq j \leq m$. The above distribution is designed to take into account the spatial structure. In that paper, they propose to use $w_{ij} = 1$ if locations $i$ and $j$ are considered neighbors. They also mention the possibility of using the inverse of the correlation matrix as the weights function. We argue that the weights chosen by the cross-validation procedure can discover the underlying spatial structure without any parametric assumption. Thus those weights might be helpful in selecting an appropriate distribution that models the underlying spatial structure. Further analysis is needed if one wants to fully compare the performances of the WLE, the MLE and the Bayesian estimator in the context of disease mapping.

## APPENDIX

Proof of Lemma 2.1. Observe that

$$\overline{X}_{i\cdot}^{(-j)} = e_n \overline{X}_{i\cdot} - \frac{1}{n-1} X_{ij},$$

where $e_n = \frac{n}{n-1}$.

Let $S_1^e = \frac{1}{n} \sum_{j=1}^{n} (\overline{X}_{1\cdot}^{(-j)} - \overline{X}_{2\cdot}^{(-j)})^2$. It then follows that

$$S_1^e = \frac{1}{n} \sum_{j=1}^{n} \left( \left( e_n \overline{X}_{1\cdot} - \frac{1}{n-1} X_{1j} \right) - \left( e_n \overline{X}_{2\cdot} - \frac{1}{n-1} X_{2j} \right) \right)^2$$

$$= \frac{1}{n}\Bigg(ne_n^2(\overline{X}_{1\cdot} - \overline{X}_{2\cdot})^2 - 2\frac{e_n}{n-1}(\overline{X}_{1\cdot} - \overline{X}_{2\cdot})\sum_{j=1}^{n}(X_{1j} - X_{2j}) + \left(\frac{1}{n-1}\right)^2\sum_{j=1}^{n}(X_{1j} - X_{2j})^2\Bigg)$$

$$= \frac{n(n-2)}{(n-1)^2}(\overline{X}_{1\cdot} - \overline{X}_{2\cdot})^2 + \frac{1}{n(n-1)^2}\sum_{j=1}^{n}(X_{1j} - X_{2j})^2.$$

Let $S_2^e = \frac{1}{n}\sum_{j=1}^{n}(\overline{X}_{1\cdot}^{(-j)} - X_{2\cdot}^{(-j)})(\overline{X}_{1\cdot}^{(-j)} - X_{1j})$. It follows that

$$S_2^e = \frac{1}{n}\sum_{j=1}^{n}\Bigg(e_n(\overline{X}_{1\cdot} - \overline{X}_{2\cdot}) - \frac{1}{n-1}(X_{1j} - X_{2j})\Bigg)\Bigg(\Bigg(e_n\overline{X}_{1\cdot} - \frac{1}{n-1}X_{1j}\Bigg) - X_{1j}\Bigg)$$

$$= \frac{1}{n}\sum_{j=1}^{n}\Bigg(e_n(\overline{X}_{1\cdot} - \overline{X}_{2\cdot}) - \frac{1}{n-1}(X_{1j} - X_{2j})\Bigg)(e_n\overline{X}_{1\cdot} - e_nX_{1j})$$

$$= \frac{e_n^2}{n}(\overline{X}_{1\cdot} - \overline{X}_{2\cdot})\sum_{j=1}^{n}(\overline{X}_{1\cdot} - X_{1j}) - \frac{e_n}{n(n-1)}\sum_{j=1}^{n}(X_{1j} - X_{2j})(\overline{X}_{1\cdot} - X_{1j})$$

$$= -\frac{e_n}{n(n-1)}\sum_{j=1}^{n}(X_{1j} - X_{2j})(\overline{X}_{1\cdot} - X_{1j}) \qquad \left[\text{since } \sum_{j=1}^{n}(\overline{X}_{1\cdot} - X_{1j}) = 0\right]$$

$$= -\frac{e_n}{n(n-1)}\Bigg(\overline{X}_{1\cdot}\sum_{j=1}^{n}(X_{1j} - X_{2j}) - \sum_{j=1}^{n}X_{1j}^2 + \sum_{j=1}^{n}X_{1j}X_{2j}\Bigg)$$

$$= -\frac{e_n}{n-1}\Bigg(\overline{X}_{1\cdot}(\overline{X}_{1\cdot} - \overline{X}_{2\cdot}) - \frac{1}{n}\sum_{j=1}^{n}X_{1j}^2 + \frac{1}{n}\sum_{j=1}^{n}X_{1j}X_{2j}\Bigg)$$

$$= \frac{n}{(n-1)^2}(\hat{\sigma}_1^2 - \widehat{\mathrm{cov}}).$$

This completes the proof. □

Proof of Proposition 2.1. By the weak law of large numbers, it follows that

$$\hat{\sigma}_1^2 - \widehat{\mathrm{cov}} \longrightarrow \sigma_1^2 - \rho\sigma_1\sigma_2.$$

Thus condition $\rho < \sigma_1/\sigma_2$ implies that $\hat{\sigma}_1^2 > \widehat{\mathrm{cov}}$ for sufficiently large $n$. Thus, $\lambda_2^{\mathrm{opt}}$ eventually will be positive. □

Proof of Proposition 2.2. From Lemma 2.1, it follows that the second term of $S_1$ goes to zero in probability as $n$ goes to infinity, while the first term converges to $(\theta_1^0 - \theta_2^0)^2$ in probability. Therefore we have

$$S_1^e \xrightarrow{P_{\boldsymbol{\theta}^0}} (\theta_1^0 - \theta_2^0)^2 \qquad \text{as } n \to \infty,$$

where $(\theta_1^0 - \theta_2^0)^2 \neq 0$ by assumption.

Moreover, we see that $S_2^e = O_P(\frac{1}{n})$. By definition of $\lambda_2^{\text{opt}}$, it follows that

$$|\lambda_2^*| = \left|\frac{S_2^e}{S_1^e}\right| \xrightarrow{P_{\boldsymbol{\theta}^0}} 0 \qquad \text{as } n \to \infty.$$

This completes the proof. $\square$

Proof of Proposition 2.3. By differentiating $D_e^{(m)} - \nu(\mathbf{1}^t\boldsymbol{\lambda} - 1)$ and setting the result to zero, it follows that

$$\frac{\partial D_e^{(m)} - \nu(\mathbf{1}^t\boldsymbol{\lambda} - 1)}{\partial \lambda} = -2b_e + 2A_e\boldsymbol{\lambda}_e^{\text{opt}} - \nu\mathbf{1} = 0.$$

It then follows that

$$\boldsymbol{\lambda}_e^{\text{opt}} = A_e^{-1}\left(b_e + \frac{\nu}{2}\mathbf{1}\right).$$

We then have

$$1 = \mathbf{1}^t\boldsymbol{\lambda}_e^{\text{opt}} = \mathbf{1}^t A_e^{-1}\left(b_e + \frac{\nu\mathbf{1}}{2}\right).$$

Thus

$$\nu = \frac{2}{\mathbf{1}^t A_e^{-1}\mathbf{1}}(1 - \mathbf{1}^t A_e^{-1} b_e).$$

Therefore

$$\boldsymbol{\lambda}_e^{\text{opt}} = A_e^{-1}\left(b_e + \frac{1 - \mathbf{1}^t A_e^{-1} b_e}{\mathbf{1}^t A_e^{-1}\mathbf{1}}\mathbf{1}\right).$$

Since $D_e^{(m)}$ is a quadratic function of $\boldsymbol{\lambda}$ and $A \geq 0$, the minimum is achieved at the point $\boldsymbol{\lambda}_e^{\text{opt}}$. Furthermore, by (5) and (7) we have

$$A_e^{-1}b_e = A_e^{-1}(A_1 - e_n^2\widehat{\Sigma}_1) = (1, 0, 0, \ldots, 0)^t - e_n^2 A_e^{-1}\widehat{\Sigma}_1.$$

Denote the optimum weight vector by $\lambda^{\text{opt}}$. It follows that

$$\boldsymbol{\lambda}_e^{\text{opt}} = (1, 0, 0, \ldots, 0)^t - e_n^2\left(A_e^{-1}\widehat{\Sigma}_1 - \frac{\mathbf{1}^t A_e^{-1}\widehat{\Sigma}_1}{\mathbf{1}^t A_e^{-1}\mathbf{1}} A_e^{-1}\mathbf{1}\right).$$

This completes the proof. $\square$

Proof of Proposition 2.4. From (8), it follows that

$$\lambda_1^{\text{opt}} = 1 - \frac{(n_1/(n_1 - 1))^2\hat{\sigma}_1^2}{n_1(\overline{X}_{1\cdot} - \overline{X}_{2\cdot})^2 + (1/(n_1 - 1))\hat{\sigma}_1^2}.$$

By the weak law of large numbers, we have

$$\hat{\sigma}_1^2 \xrightarrow{P_{\boldsymbol{\theta}^0}} \sigma_1^2,$$

$$(\overline{X}_{1\cdot} - \overline{X}_{2\cdot})^2 \xrightarrow{P_{\boldsymbol{\theta}^0}} (\theta_1^0 - \theta_2^0)^2 \neq 0.$$

It then follows that

$$\frac{(n_1/(n_1-1))^2\hat{\sigma}_1^2}{n_1(\overline{X}_{1\cdot} - \overline{X}_{2\cdot})^2 + (1/(n_1-1))\hat{\sigma}_1^2} \xrightarrow{P_{\boldsymbol{\theta}^0}} 0.$$

We then have

$$\lambda_1^{\mathrm{opt}} \xrightarrow{P_{\boldsymbol{\theta}^0}} 1.$$

The last assertion of the theorem follows by the fact that $\lambda_1 + \lambda_2 = 1$. □

Proof of Theorem 3.1. Consider

$$\begin{aligned}
\frac{1}{n_1}D_{n_1}(\boldsymbol{\lambda}) &= \frac{1}{n_1}\sum_{j=1}^{n_1}(X_{1j} - \phi(\tilde{\theta}_1^{(-j)}))^2 \\
&= \frac{1}{n_1}\sum_{j=1}^{n_1}((X_{1j} - \phi(\hat{\theta}_1^{(-j)})) + (\phi(\hat{\theta}_1^{(-j)}) - \phi(\tilde{\theta}_1^{(-j)})))^2 \\
&= \frac{1}{n_1}\sum_{j=1}^{n_1}(X_{1j} - \phi(\hat{\theta}_1^{(-j)}))^2 + \frac{1}{n_1}\sum_{j=1}^{n_1}(\phi(\hat{\theta}_1^{(-j)}) - \phi(\tilde{\theta}_1^{(-j)}))^2 \\
&\quad + \frac{2}{n_1}\sum_{j=1}^{n_1}(X_{1j} - \phi(\hat{\theta}_1^{(-j)}))(\phi(\hat{\theta}_1^{(-j)}) - \phi(\tilde{\theta}_1^{(-j)})).
\end{aligned}$$

Note that

$$\begin{aligned}
&\frac{1}{n_1}\sum_{j=1}^{n_1}(X_{1j} - \phi(\hat{\theta}_1^{(-j)}))(\phi(\hat{\theta}_1^{(-j)}) - \phi(\tilde{\theta}_1^{(-j)})) \\
&\quad = \frac{1}{n_1}\sum_{j=1}^{n_1}(X_{1j} - \phi(\theta_1^0))(\phi(\hat{\theta}_1^{(-j)}) - \phi(\tilde{\theta}_1^{(-j)})) \\
&\qquad + \frac{1}{n_1}\sum_{j=1}^{n_1}(\phi(\theta_1^0) - \phi(\hat{\theta}_1^{(-j)}))(\phi(\hat{\theta}_1^{(-j)}) - \phi(\tilde{\theta}_1^{(-j)})) \\
&\quad = S_1 + S_2,
\end{aligned}$$

where

$$S_1 = \frac{1}{n_1}\sum_{j=1}^{n_1}(X_{1j} - \phi(\theta_1^0))(\phi(\hat{\theta}_1^{(-j)}) - \phi(\tilde{\theta}_1^{(-j)})),$$

$$S_2 = \frac{1}{n_1}\sum_{j=1}^{n_1}(\phi(\theta_1^0) - \phi(\hat{\theta}_1^{(-j)}))(\phi(\hat{\theta}_1^{(-j)}) - \phi(\tilde{\theta}_1^{(-j)})).$$

We first show that $S_1 \xrightarrow{P_{\boldsymbol{\theta}^0}} 0$.

Consider

$$\begin{aligned}
&P_{\boldsymbol{\theta}^0}(|S_1| > \varepsilon) \\
&\quad = P_{\boldsymbol{\theta}^0}(\varepsilon < |S_1| \text{ and } |\phi(\hat{\theta}_1^{(-j)}) - \phi(\tilde{\theta}_1^{(-j)})| < M \text{ for all } j) \\
&\qquad + P_{\boldsymbol{\theta}^0}(\varepsilon < |S_1| \text{ and } |\phi(\hat{\theta}_1^{(-l)}) - \phi(\tilde{\theta}_1^{(-l)})| \ge M \text{ for some } l) \\
&\quad \le P_{\boldsymbol{\theta}^0}\left(\varepsilon < |S_1| < \frac{M}{n_1}\sum_{j=1}^{n_1}|X_{1j} - \phi(\theta_1^0)|\right) \\
&\qquad + \sum_{l=1}^{n_1} P_{\boldsymbol{\theta}^0}(|\phi(\hat{\theta}_1^{(-l)}) - \phi(\tilde{\theta}_1^{(-l)})| \ge M) \\
&\quad \le P_{\boldsymbol{\theta}^0}\left(\frac{\varepsilon}{M} < \left|\frac{1}{n_1}\sum_{j=1}^{n_1}(X_{1j} - \phi(\theta_1^0))\right|\right) + n_1 P_{\boldsymbol{\theta}^0}(|\phi(\hat{\theta}_1^{(-1)}) - \phi(\tilde{\theta}_1^{(-1)})| \ge M) \\
&\quad = P_{\boldsymbol{\theta}^0}\left(\left|\frac{1}{n_1}\sum_{j=1}^{n_1}(X_{1j} - \phi(\theta_1^0))\right| > \frac{1}{M}\varepsilon\right) \\
&\qquad + n_1 P_{\boldsymbol{\theta}^0}(|\phi(\hat{\theta}_1^{(n_1-1)}) - \phi(\tilde{\theta}_1^{(n_1-1)})| \ge M).
\end{aligned}$$

The first term goes to zero by the weak law of large numbers. The second term also goes to zero by assumption (iv). We then have

$$P_{\boldsymbol{\theta}^0}(|S_1| > \varepsilon) \longrightarrow 0 \qquad \text{as } n_1 \to \infty. \tag{18}$$

We next show that $S_2 \xrightarrow{P_{\boldsymbol{\theta}^0}} 0$ as $n_1 \to \infty$.

Consider

$$\begin{aligned}
&P_{\boldsymbol{\theta}^0}(|S_2| > \varepsilon) \\
&\quad = P_{\boldsymbol{\theta}^0}(\varepsilon < |S_2| \text{ and } |\phi(\hat{\theta}_1^{(-j)}) - \phi(\tilde{\theta}_1^{(-j)})| < M \text{ for all } j) \\
&\qquad + P_{\boldsymbol{\theta}^0}(\varepsilon < |S_2| \text{ and } |\phi(\hat{\theta}_1^{(-l)}) - \phi(\tilde{\theta}_1^{(-l)})| \ge M \text{ for some } l) \\
&\quad \le P_{\boldsymbol{\theta}^0}\left(\varepsilon < |S_2| < \frac{M}{n_1}\left|\sum_{j=1}^{n_1}(\phi(\hat{\theta}_1^{(-j)}) - \phi(\theta_1^0))\right|\right) \\
&\qquad + \sum_{l=1}^{n_1} P_{\boldsymbol{\theta}^0}(|\phi(\hat{\theta}_1^{(-l)}) - \phi(\tilde{\theta}_1^{(-l)})| \ge M)
\end{aligned}$$

$$
\begin{aligned}
&\leq P_{\boldsymbol{\theta}^0}\left(\frac{1}{M}\varepsilon < \left|\frac{1}{n_1}\sum_{j=1}^{n_1}(\phi(\hat{\theta}_1^{(-j)}) - \phi(\theta_1^0))\right|\right) \\
&\quad + n_1 P_{\boldsymbol{\theta}^0}(|\phi(\hat{\theta}_1^{(-1)}) - \phi(\tilde{\theta}_1^{(-1)})| \geq M) \\
&= P_{\boldsymbol{\theta}^0}\left(\left|\frac{1}{n_1}\sum_{j=1}^{n_1}(\phi(\hat{\theta}_1^{(-j)}) - \phi(\theta_1^0))\right| > \frac{1}{M}\varepsilon\right) \\
&\quad + n_1 P_{\boldsymbol{\theta}^0}(|\phi(\hat{\theta}_1^{(n_1-1)}) - \phi(\tilde{\theta}_1^{(n_1-1)})| \geq M).
\end{aligned}
$$

The first term goes to zero by assumption (ii). The second term also goes to zero by assumption (iv). We then have

$$
P_{\boldsymbol{\theta}^0}(|S_2| > \varepsilon) \longrightarrow 0 \qquad \text{as } n_1 \to \infty. \tag{19}
$$

It then follows that

$$
\begin{aligned}
\frac{1}{n_1}D_{n_1}(\boldsymbol{\lambda}) &= \frac{1}{n_1}\sum_{j=1}^{n_1}(X_{1j} - \phi(\hat{\theta}_1^{(-j)}))^2 \\
&\quad + \frac{1}{n_1}\sum_{j=1}^{n_1}(\phi(\hat{\theta}_1^{(-j)}) - \phi(\tilde{\theta}_1^{(-j)}))^2 + R_n,
\end{aligned} \tag{20}
$$

where $R_n \xrightarrow{P_{\boldsymbol{\theta}^0}} 0$. Observe that the first term is independent of $\boldsymbol{\lambda}$. Therefore the second term must be minimized with respect to $\boldsymbol{\lambda}$ to obtain the minimum of $\frac{1}{n_1}D_{n_1}(\boldsymbol{\lambda})$. We see that the second term is always nonnegative. It then follows that, with probability tending to 1,

$$
\frac{1}{n_1}D_{n_1}(\boldsymbol{\lambda}) \geq \frac{1}{n_1}\sum_{j=1}^{n_1}(X_{1j} - \phi(\hat{\theta}_1^{(-j)}))^2 = \frac{1}{n_1}D_{n_1}(\mathbf{w}),
$$

since $\phi(\hat{\theta}_1^{(-j)}) = \phi(\tilde{\theta}_1^{(-j)})$ for $\boldsymbol{\lambda}^{(cv)} = \mathbf{w}_0 = (1,0,0,\dots,0)^t$ for fixed $n_1$.

Finally, we will show that

$$
\boldsymbol{\lambda}^{(cv)} \xrightarrow{P_{\boldsymbol{\theta}^0}} \mathbf{w}_0 \qquad \text{as } n_1 \to \infty.
$$

Suppose to the contrary that $\boldsymbol{\lambda}^{(cv)} \xrightarrow{P_{\boldsymbol{\theta}^0}} \mathbf{w}_0 + \mathbf{d}$ where $\mathbf{d}$ is a nonzero vector. Then there exists $n_0$ such that for $n_1 > n_0$,

$$
\frac{1}{n_1}D_{n_1}(\boldsymbol{\lambda}^{(cv)}) \geq \frac{1}{n_1}D_{n_1}(\mathbf{w}).
$$

This is a contradiction because $\boldsymbol{\lambda}^{(cv)}$ is the vector which minimizes $\frac{1}{n_1}D_{n_1}$ for any fixed $n_1$ and the minimum of $\frac{1}{n_1}D_{n_1}(\boldsymbol{\lambda})$ is unique by assumption. □

Proof of Lemma 3.1. Recall that the average discrepancy of cross-validation for the log-normal case is given by

$$\frac{1}{n}D_n(\lambda_1,\lambda_2) = \frac{1}{n}\sum_{j=1}^{n}(X_{1j} - e^{(\lambda_1/(n-1))\sum_{k\neq j}\log(x_{1k})+(\lambda_2/(n-1))\sum_{k\neq j}\log(x_{2k})+1/2})^2. \tag{21}$$

(i) Since we require that $\lambda_1+\lambda_2=1$, we can rewrite the average discrepancy as

$$\frac{1}{n}D_n(1-\lambda_2,\lambda_2) = \frac{1}{n}\sum_{j=1}^{n}(X_{1j} - e^{\overline{Y}_{1\cdot}^{(-j)}+\lambda_2(\overline{Y}_{2\cdot}^{(-j)}-\overline{Y}_{1\cdot}^{(-j)})+1/2})^2, \tag{22}$$

where $Y_{ij}=\log(X_{ij})$, $i=1,2; j=1,2,\dots,n$.

Note that $\alpha(x)=(x-a)^2$ and $\beta(x)=e^{b*x+c}$ are both convex functions for any given constants $a,b$ and $c$. It then follows that $\gamma(x)=(e^{b*x+c}-a)^2$ is also a convex function. Thus, $\frac{1}{n}D_n(1-\lambda_2,\lambda_2)$ is a strict convex function with respect to $\lambda_2$ for fixed $n$.

(ii) The first-order derivative of $\frac{1}{n}D_n$ is given by

$$\begin{aligned}\frac{1}{n}\frac{\partial D_n(1-\lambda_2,\lambda_2)}{\partial\lambda_2} &= -\frac{2}{n}\sum_{j=1}^{n}(X_{1j} - e^{\overline{Y}_{1\cdot}^{(-j)}+\lambda_2(\overline{Y}_{2\cdot}^{(-j)}-\overline{Y}_{1\cdot}^{(-j)})+1/2}) \\ &\qquad * e^{\overline{Y}_{1\cdot}^{(-j)}+\lambda_2(\overline{Y}_{2\cdot}^{(-j)}-\overline{Y}_{1\cdot}^{(-j)})+1/2} * (\overline{Y}_{2\cdot}^{(-j)}-\overline{Y}_{1\cdot}^{(-j)}).\end{aligned} \tag{23}$$

Observe that

$$\overline{Y}_{2\cdot}^{(-j)}-\overline{Y}_{1\cdot}^{(-j)} = (\overline{Y}_{2\cdot}-\overline{Y}_{1\cdot}) + \frac{1}{n-1}([\overline{Y}_{2\cdot}-\overline{Y}_{1\cdot}]-[Y_{2j}-Y_{1j}]).$$

It then follows that

$$\begin{aligned}\frac{1}{n}\frac{\partial D_n(1-\lambda_2,\lambda_2)}{\partial\lambda_2} &= -\frac{2}{n}\sum_{j=1}^{n}(X_{1j} - e^{\overline{Y}_{1\cdot}+\lambda_2(\overline{Y}_{2\cdot}-\overline{Y}_{1\cdot})+T_j^n+1/2}) \\ &\qquad * e^{\overline{Y}_{1\cdot}+\lambda_2(\overline{Y}_{2\cdot}-\overline{Y}_{1\cdot})+T_j^n+1/2} * ((\overline{Y}_{2\cdot}-\overline{Y}_{1\cdot})+R_j^n),\end{aligned} \tag{24}$$

where

$$\begin{aligned}R_j^n &= \frac{1}{n-1}([\overline{Y}_{2\cdot}-\overline{Y}_{1\cdot}]-[Y_{2j}-Y_{1j}]) = O_P(n^{-1}), \\ T_j^n(\lambda_2) &= \lambda_2 R_j^n + \frac{1}{n-1}(\overline{Y}_{1\cdot}-Y_{1j}).\end{aligned} \tag{25}$$

Thus

$$\frac{1}{n}\frac{\partial D_n(1-\lambda_2,\lambda_2)}{\partial\lambda_2} = -2F_n(\lambda_2) * E_n(\lambda_2), \tag{26}$$

where

$$F_n(\lambda_2) = e^{\overline{Y}_{1\cdot}+\lambda_2(\overline{Y}_{2\cdot}-\overline{Y}_{1\cdot})+1/2} * (\overline{Y}_{2\cdot}-\overline{Y}_{1\cdot}) \tag{27}$$

and

$$\begin{aligned} E_n(\lambda_2) = \frac{1}{n}\sum_{j=1}^{n}&(X_{1j} - e^{\overline{Y}_{1\cdot}+\lambda_2(\overline{Y}_{2\cdot}-\overline{Y}_{1\cdot})+T_j^n+1/2}) \\ &* (1+R_j^n/(\overline{Y}_{2\cdot}-\overline{Y}_{1\cdot})) * e^{T_j^n}. \end{aligned} \tag{28}$$

For any $|\lambda_2| \le 1$, we have $T_j^n(\lambda_2) = O_P(n^{-1})$ and $e^{T_j^n} = 1 + T_j^n + O_P(n^{-2})$. Thus

$$\begin{aligned} E_n(\lambda_2) = \frac{1}{n}\sum_{j=1}^{n}&(X_{1j} - e^{\overline{Y}_{1\cdot}+\lambda_2(\overline{Y}_{2\cdot}-\overline{Y}_{1\cdot})+T_j^n+1/2}) \\ &* (1+R_j^n/(\overline{Y}_{1\cdot}-\overline{Y}_{2\cdot})) \\ &* (1+T_j^n+O_P(n^{-2})), \qquad |\lambda_2| \le 1. \end{aligned} \tag{29}$$

Furthermore, for any $|\lambda_2| \le 1$ we have

$$\begin{aligned} E_n(\lambda_2) = \frac{1}{n}\sum_{j=1}^{n}&(X_{1j} - e^{\overline{Y}_{1\cdot}+\lambda_2(\overline{Y}_{2\cdot}-\overline{Y}_{1\cdot})+T_j^n+1/2}) \\ &+ U_n(\lambda_2)/(\overline{Y}_{2\cdot}-\overline{Y}_{1\cdot}) + V_n(\lambda_2) + W_n(\lambda_2)/(\overline{Y}_{2\cdot}-\overline{Y}_{1\cdot}), \end{aligned} \tag{30}$$

where

$$U_n(\lambda_2) = \frac{1}{n}\sum_{j=1}^{n}(X_{1j} - e^{\overline{Y}_{1\cdot}^{(-j)}+\lambda_2(\overline{Y}_{2\cdot}^{(-j)}-\overline{Y}_{1\cdot}^{(-j)})+1/2}) * R_j^n,$$

$$V_n(\lambda_2) = \frac{1}{n}\sum_{j=1}^{n}(X_{1j} - e^{\overline{Y}_{1\cdot}^{(-j)}+\lambda_2(\overline{Y}_{2\cdot}^{(-j)}-\overline{Y}_{1\cdot}^{(-j)})+1/2}) * (T_j^n + O_P(n^{-2})),$$

$$W_n(\lambda_2) = \frac{1}{n}\sum_{j=1}^{n}(X_{1j} - e^{\overline{Y}_{1\cdot}^{(-j)}+\lambda_2(\overline{Y}_{2\cdot}^{(-j)}-\overline{Y}_{1\cdot}^{(-j)})+1/2}) * (T_j^n + O_P(n^{-2})) * R_j^n.$$

If $|\lambda_2| < 1$, then

$$|\overline{Y}_{1\cdot}^{(-j)} + \lambda_2(\overline{Y}_{2\cdot}^{(-j)} - \overline{Y}_{1\cdot}^{(-j)})| \le |\overline{Y}_{1\cdot}^{(-j)}| + |\overline{Y}_{2\cdot}^{(-j)} - \overline{Y}_{1\cdot}^{(-j)}|. \tag{31}$$

We also consider

$$\begin{aligned} B_n(\lambda_2) &= \frac{1}{n}\sum_{j=1}^{n} T_j^n(\lambda_2) \\ &= \frac{1}{n}\sum_{j=1}^{n}\left(\frac{1}{n-1}\left[\lambda_2((\overline{Y}_{1\cdot}-\overline{Y}_{2\cdot}) - (Y_{1j}-Y_{2j})) - \frac{1}{n-1}(\overline{Y}_{1\cdot}-Y_{1j})\right]\right). \end{aligned}$$

Note that, for any sequence of random variables $Z_j$, $j=1,2,\dots,n$, $E|(Z_iZ_k)|<\infty$, $i,k=1,2,\dots,n$,

$$P\left(\left|\frac{1}{n(n-1)}\sum_{j=1}^{n}Z_n\right|>\varepsilon\right)\le\frac{1}{n^4}E\left(\sum_i\sum_k Z_iZ_k\right)=O(n^{-2}). \tag{32}$$

By combining (25), (31) and (32), we can show that, for $|\lambda_2|\le 1$,

$$\begin{aligned} U_n(\lambda_2)&=O_P(n^{-2}), & V_n(\lambda_2)&=O_P(n^{-2}),\\ W_n(\lambda_2)&=O_P(n^{-2}), & B_n(\lambda_2)&=O_P(n^{-2}). \end{aligned} \tag{33}$$

We also observe that

$$\begin{aligned} &\frac{1}{n}\sum_{j=1}^{n}(X_{1j}-e^{\overline{Y}_{1\cdot}+\lambda_2(\overline{Y}_{2\cdot}-\overline{Y}_{1\cdot})+T_j^n+1/2})\\ &\qquad=\frac{1}{n}\sum_{j=1}^{n}(X_{1j}-e^{\overline{Y}_{1\cdot}+\lambda_2(\overline{Y}_{2\cdot}-\overline{Y}_{1\cdot})+1/2}(1+T_j^n(\lambda_2)+O_P(n^{-2}))). \end{aligned}$$

It then follows that

$$E_n(\lambda_2)=\frac{1}{n}\sum_{j=1}^{n}(X_{1j}-e^{\overline{Y}_{1\cdot}+\lambda_2(\overline{Y}_{2\cdot}-\overline{Y}_{1\cdot})+1/2})+C_n(\lambda_2), \tag{34}$$

where $C_n=B_n(\lambda_2)+U_n(\lambda_2)/(\overline{Y}_{2\cdot}-\overline{Y}_{1\cdot})+V_n(\lambda_2)+W_n(\lambda_2)/(\overline{Y}_{2\cdot}-\overline{Y}_{1\cdot})$.

It is clear that $\overline{Y}_{2\cdot}-\overline{Y}_{1\cdot}\xrightarrow{\text{a.s.}}\mu_1^0-\mu_2^0$. Thus

$$E_n(\lambda_2)=\frac{1}{n}\sum_{j=1}^{n}X_{1j}-e^{\overline{Y}_{1\cdot}+\lambda_2(\overline{Y}_{2\cdot}-\overline{Y}_{1\cdot})+1/2}+O_P(n^{-2}). \tag{35}$$

Without essential loss of generality, we assume that $\mu_1^0>\mu_2^0$. It then follows that

$$E_n(1)\xrightarrow{\text{a.s.}}e^{\mu_1^0+1/2}(1-e^{\mu_2^0-\mu_1^0})<0 \tag{36}$$

and

$$E_n(-1)\xrightarrow{\text{a.s.}}e^{\mu_1^0+1/2}(1-e^{\mu_1^0-\mu_2^0})>0. \tag{37}$$

By (26), (27), (36) and (37), it follows that for sufficiently large $n$,

$$\begin{aligned} &\frac{4}{n^2}*\frac{\partial D_n(1-\lambda_2,\lambda_2)}{\partial\lambda_2}\bigg|_{\lambda_2=1}*\frac{\partial D_n(1-\lambda_2,\lambda_2)}{\partial\lambda_2}\bigg|_{\lambda_2=-1}\\ &\quad=F_n(1)*F_n(-1)*E_n(1)*E(-1)<0. \end{aligned} \tag{38}$$

Since $D_n$ is strictly convex, then its second-order derivative is positive. Therefore, the first-order derivative of $D_n$ is monotone. By (38), we then have that the optimal weight $\lambda_2^*\in(-1,1)$ for sufficiently large $n$ with probability tending to 1. Furthermore, it converges to a unique limit. Suppose that this is not true and there are two limits $\overline{\lambda}_2^I$ and $\overline{\lambda}_2^{II}$. Then $0.5\overline{\lambda}_2^I + 0.5\overline{\lambda}_2^{II}$ achieves a small value for $\frac{1}{n}D_n(1-\lambda_2,\lambda_2)$ since it is strictly convex. This is a contradiction. $\square$

Proof of Corollary 3.2. It suffices to show the assumptions of Theorem 3.1 are satisfied for the log-normal case.

(i) By Lemma 3.1, $\frac{1}{n}D_n(1-\lambda_2,\lambda_2)$ is a strict convex function with respect to $\lambda_2$. Therefore assumption (i) of Theorem 3.1 is satisfied.

(ii) We then check assumption (ii) of Theorem 3.1. Let $\frac{1}{n}S_n^I = \frac{1}{n}\times\sum_{j=1}^n(\phi(\mu_1^{(-j)}) - \phi(\mu_1^0))$. Thus

$$\frac{1}{n}S_n^I = \frac{1}{n}\sum_{j=1}^n (e^{(1/(n-1))\sum_{k\neq j}\log(X_{1k})+1/2} - e^{\mu_1^0+1/2}).$$

Let $A_j^n = e^{(1/(n-1))\sum_{k\neq j}\log(X_{1k})+1/2} - e^{\mu_1^0+1/2}$. It then follows that

$$\frac{1}{n}S_n^I = \frac{1}{n}\sum_{j=1}^n A_j^n.$$

Observe that $Y_{ij} = \log(X_{ij}) \sim N(\mu_i^0, 1)$, $j = 1, 2, \ldots, n$. Thus we have

$$E_{\mu_1^0}(e^{\log(X_{ij})*t}) = E(e^{Y_{ij}*t}) = e^{\mu_i^0 t + t^2/2}, \tag{39}$$

for $i = 1, 2; j = 1, 2, \ldots, n$. We then have

$$\begin{aligned} E_{\mu_1^0} e^{(1/(n-1))\sum_{k=2}^n \log(X_{1k})} &= (e^{1/(n-1)\mu_1^0 + 1/(2(n-1)^2)})^{n-1} \\ &= e^{\mu_1^0 + 1/(2(n-1))}. \end{aligned} \tag{40}$$

We also have

$$\begin{aligned} &E_{\mu_1^0} e^{(1/(n-1))\sum_{k=2}^n \log(X_{1k})} * e^{(1/(n-1))\sum_{l=1}^{n-1}\log(X_{1l})} \\ &\quad = E_{\mu_1^0}(e^{(1/(n-1))\log(X_{11}) + (1/(n-1))\log(X_{1n})}) * E(e^{(2/(n-2))\sum_{k=2}^{n-1}\log(X_{1k})}) \\ &\quad = e^{2*(1/(n-1))\mu_1^0 + 1/(2(n-1)^2)} * e^{(n-2)*(2/(n-1))\mu_1^0 + 2/(n-1)^2)} \\ &\quad = e^{2\mu_1^0 + ((2n-4)/(n-1)^2)}. \end{aligned} \tag{41}$$

By (39) and (40), it then follows that

$$\begin{aligned} E_{\mu_1^0}(A_1^n)^2 &= e * E_{\mu_1^0}[e^{(1/(n-1))\sum_{k=2}^n \log(X_{1k})} - e^{\mu_1^0}]^2 \\ &= e * [E_{\mu_1^0}(e^{(1/(n-1))\sum_{k=2}^n \log(X_{1k})})^2 \\ &\qquad - 2e^{\mu_1^0} * E_{\mu_1^0}(e^{(1/(n-1))\sum_{k=2}^n \log(X_{1k})}) + e^{2\mu_1^0}] \\ &= e * [E_{\mu_1^0}(e^{(2/n-1))\sum_{k=2}^n \log(X_{1k})}) - 2e^{\mu_1^0} * e^{\mu_1^0 + 1/(2(n-1))} + e^{2\mu_1^0}] \end{aligned}$$

$$= e * [(E_{\mu_1^0} e^{(2/(n-1))\log(X_{11})})^{n-1} - 2e^{\mu_1^0} * e^{\mu_1^0+1/(2(n-1))} + e^{2\mu_1^0}]$$

$$= e * [e^{2\mu_1^0+2/(n-1)} - 2e^{\mu_1^0} * e^{\mu_1^0+1/(2(n-1))} + e^{2\mu_1^0}] \qquad \text{[by (15)]}$$

$$= e^{2\mu_1^0+1} O\left(\frac{1}{n}\right).$$

By (40) and (41), we also have

$$\begin{aligned} E_{\mu_1^0}(A_1^n * A_n^n) &= E_{\mu_1^0}(e^{(1/(n-1))\sum_{k=2}^n \log(X_{1k})+1/2} - e^{\mu_1^0+1/2}) \\ &\quad * (e^{(1/(n-1))\sum_{l=1}^{n-1} \log(X_{1k})+1/2} - e^{\mu_1^0+1/2}) \\ &= e(E_{\mu_1^0} e^{(1/(n-1))\sum_{k=2}^n \log(X_{1k})} * e^{(1/(n-1))\sum_{l=1}^{n-1} \log(X_{1l})} \\ &\qquad - 2e^{\mu_1^0} * E_{\mu_1^0} e^{(1/(n-1))\sum_{k=2}^n \log(X_{1k})} + e^{2\mu_1^0}) \\ &= e * (e^{2\mu_1^0+((2n-4)/(n-1)^2)} - 2 * e^{\mu_1^0} * e^{\mu_1^0+1/(2(n-1))} + e^{2\mu_1^0}) \\ &= e^{2\mu_1^0+1}(e^{((2n-4)/(n-1)^2)} - 2 * e^{1/(2(n-1))} + 1) \\ &= e^{2\mu_1^0+1} O\left(\frac{1}{n}\right). \end{aligned}$$

For any fixed $j$ and $k$, we then have

$$E_{\mu_1^0}(A_j^n * A_k^n) = O\left(\frac{1}{n}\right). \tag{42}$$

Therefore,

$$\begin{aligned} P_{\mu_1^0}\left(\left|\frac{1}{n} S_n^I\right| > \varepsilon\right) &\le \frac{1}{n^2\varepsilon^2} E(S_n^I)^2 \\ &= \frac{1}{n^2\varepsilon^2} E\left(\sum_{j=1}^n (A_j^n)^2\right) + \frac{1}{n^2\varepsilon^2}\left(\sum_{j=1}^n \sum_{k\neq j} E_{\mu_1^0}(A_j^n * A_k^n)\right) \\ &= \frac{1}{n\varepsilon^2} E(A_1^n)^2 + \frac{n(n-1)}{n^2\varepsilon^2} E_{\mu_1^0}(A_1^n A_n^n) \\ &= O\left(\frac{1}{n}\right) \longrightarrow 0 \qquad \text{as } n \to \infty. \end{aligned}$$

This implies that assumption (ii) is satisfied for the log-normal case.

(iii) Let

$$\frac{1}{n} S_n^{II} = \frac{1}{n} \sum_{j=1}^n (X_{1j} - \phi(\hat{\mu}_1^{(-j)})^2).$$

Observe that

$$\begin{aligned}\frac{1}{n}S_n^{II} &= \frac{1}{n}\sum_{j=1}^{n} X_{1j}^2 - \frac{1}{n}\sum_{j=1}^{n} X_{1j}\phi(\hat{\mu}_1^{(-j)}) + \frac{1}{n}\sum_{j=1}^{n}\phi^2(\hat{\mu}_1^{(-j)})\\ &= I_1^n + I_2^n + I_3^n,\end{aligned}$$

where $I_1^n = \frac{1}{n}\sum_{j=1}^n X_{1j}^2$, $I_2^n = \frac{1}{n}\sum_{j=1}^n X_{1j}\phi(\hat{\mu}_1^{(-j)})$ and $I_3^n = \frac{1}{n}\sum_{j=1}^n \phi^2(\hat{\mu}_1^{(-j)})$. By the weak law of large numbers, it follows that

$$I_1^n = \frac{1}{n}\sum_{j=1}^{n} X_{1j}^2 \xrightarrow{P_{\mu_1^0}} E(X_{11})^2 = e^{2\mu_1^0+2}. \tag{43}$$

Consider

$$I_2^n = \frac{1}{n}\sum_{j=1}^{n} X_{1j}\phi(\hat{\mu}_1^{(-j)}) = \frac{1}{n}\sum_{j=1}^{n} e^{Y_{1j}+(1/(n-1))\sum_{k\neq j} Y_{1k}+1/2}, \tag{44}$$

where $Y_{1j} = \log(X_{1j}) \sim N(\mu_1^0, 1), j = 1, 2, \dots, n$.

Note that for any $j$

$$Y_{1j} + \frac{1}{n-1}\sum_{k\neq j} Y_{1k} = \frac{n-2}{n-1}Y_{1j} + \frac{n}{n-1}\overline{Y}_{1\cdot},$$

where $\overline{Y}_{1\cdot} = \frac{1}{n}\sum_{k=1}^n Y_{1k}$.

It then follows that

$$I_2^n = e^{1/2} * e^{(n/(n-1))\overline{Y}_{1\cdot}} * \left(\frac{1}{n}\sum_{k=1}^{n} e^{((n-2)/(n(n-1)))*Y_{1k}}\right). \tag{45}$$

Note that

$$\begin{aligned}&\frac{1}{n}\sum_{k=1}^{n} e^{((n-2)/(n(n-1)))*Y_{1k}}\\ &\qquad = \frac{1}{n}\sum_{k=1}^{n}\left(1 + \frac{n-2}{n(n-1)} * Y_{1k} + O_P(n^{-2})\right) \xrightarrow{P_{\mu_1^0}} 1.\end{aligned} \tag{46}$$

It then follows that

$$I_2^n \xrightarrow{P_{\mu_1^0}} e^{\mu_1^0+1/2} \qquad \text{as } n\to\infty. \tag{47}$$

We also have

$$I_3^n = e * e^{(2n/(n-1))\overline{Y}_{1\cdot}} \frac{1}{n}\sum_{j=1}^{n} e^{-(2/(n-2))Y_j} \xrightarrow{P_{\mu_1^0}} e^{2\mu_1^0+1}.$$

It then follows that

$$\frac{1}{n}S_2^{II} \xrightarrow{P_{\mu_1^0}} e^{2\mu_1^0+2} - 2 * e^{\mu_1^0+1/2} + e^{2\mu_1^0+1}. \tag{48}$$

This implies that assumption (iii) is satisfied.

(iv) We are now in a position to verify Assumption (iv) of Theorem 3.1. Note that the optimum weight $\lambda_2^*$ is chosen such that

$$\left.\frac{\partial D_n(1-\lambda_2,\lambda_2)}{\partial \lambda_2}\right|_{\lambda_2=\lambda_2^*}=0.$$

By (27), we see that either $F_n>0$ or $F_n<0$ for sufficiently large $n$ if $\mu_1 \neq \mu_2$. By (26), (35) and Lemma 3.1, it follows that the optimum weight $\lambda_2^*(n)$ satisfies

$$0=E_n(\lambda_2^*)=\frac{1}{n}\sum_{j=1}^{n}(X_{1j}-e^{\overline{Y}_{1\cdot}+\lambda_2^*(\overline{Y}_{2\cdot}-\overline{Y}_{1\cdot})+1/2})+O_P(n^{-2}). \tag{49}$$

We then have

$$\phi(\tilde{\mu}_1^n)=e^{\overline{Y}_{1\cdot}+\lambda_2^*(\overline{Y}_{2\cdot}-\overline{Y}_{1\cdot})+1/2}=\frac{1}{n}\sum_{j=1}^{n}X_{1j}+O_P(n^{-2}). \tag{50}$$

For sufficiently large $n$ and any constant $M>0$, say 1, and a certain $C(M)$, which depends on $M$ and whose value is of no relevance to the argument, we have

$$\begin{aligned}
&P_{\mu_1^0}(|\phi(\hat{\mu}_1^n)-\phi(\tilde{\mu}_1^n)|>M)\\
&\quad=P_{\mu_1^0}\left(\left|e^{1/n\sum_{j=1}^n\log(X_{1j})+1/2}-\frac{1}{n}\sum_{j=1}^{n}X_{1j}+O_P(n^{-2})\right|>M\right)\\
&\quad\le P_{\mu_1^0}(|e^{1/n\sum_{j=1}^n\log(X_{1j})+1/2}-e^{\mu_1^0+1/2}|>M/2)\\
&\qquad+P_{\mu_1^0}\left(\left|\frac{1}{n}\sum_{j=1}^{n}X_{1j}-e^{\mu_1^0+1/2}\right|>M/2\right)+O(n^{-2})\\
&\quad\le P_{\mu_1^0}\left(\left|\frac{1}{n}\sum_{j=1}^{n}\log(X_{1j})-\mu_1^0\right|>C(M)\right)\\
&\qquad+P_{\mu_1^0}\left(\left|\frac{1}{n}\sum_{j=1}^{n}X_{1j}-e^{\mu_1^0+1/2}\right|>M/2\right)+O(n^{-2})\\
&\quad=O(n^{-2}).
\end{aligned}$$

The last inequality follows since the fourth moments of $X_{1j}$ and $\log(X_{1j})$ both exist for any fixed $j$. Therefore, the last assumption of Theorem 3.1 is satisfied for the log-normal case. This completes the proof. □

The proofs of Theorems 3.2–3.4 resemble the proofs for fixed weights as given by Wang, van Eeden and Zidek (2004). These theorems can be proved by using Theorem 3.1 and replacing fixed weights with adaptive weights in weighted likelihood estimation. Details can be found in Wang (2001).

**Acknowledgments.** We would like to thank an Associate Editor and two anonymous referees for constructive and incisive comments on the preliminary version of this article.

Department of Mathematics
and Statistics
York University
Toronto, Ontario
Canada M3J 1P3
e-mail: stevenw@mathstat.yorku.ca

Department of Statistics
University of British Columbia
Vancouver, British Columbia
Canada V6T 1Z2
e-mail: jim@stat.ubc.ca